\documentclass[9pt]{article}
\usepackage[utf8]{inputenc}
\usepackage{amsmath}
\usepackage[numbers]{natbib} 
\usepackage{textcomp}
\usepackage{gensymb}
\usepackage{amssymb}
\usepackage{mathtools}
\usepackage{graphicx}
\usepackage{epsfig}
\usepackage{float}
\usepackage{comment}
\usepackage{caption}
\usepackage{subcaption}
\usepackage[backref=page]{hyperref}
\usepackage{cleveref}
\usepackage{array}
\usepackage{siunitx}
\usepackage{amsthm}
\usepackage{soul}
\usepackage[abs]{overpic}
\usepackage[dvipsnames]{xcolor}
\usepackage[normalem]{ulem}
\usepackage[margin=1.5in]{geometry}

\DeclareMathOperator{\sech}{sech}

\theoremstyle{definition}
\newtheorem{defn}{\bf Definition}[section]
\newtheorem{rmk}{\bf Remark}[section]

\newcommand{\sw}{}
\newcommand{\swr}{}
\newcommand{\eos}{}
\newcommand{\swb}{}
\usepackage[normalem]{ulem}
\usepackage{comment}
\begin{document}

\title{Rate-Induced Tipping to Metastable Zombie Fires}

\author{
Eoin O'Sullivan\footnote{eoin.geoffrey.osullivan@umail.ucc.ie}~\footnote{University College Cork, School of Mathematical Sciences,
Western Road, Cork T12 XF62, Ireland},  Kieran Mulchrone$^{\dagger}$ and Sebastian Wieczorek$^{\dagger}$}

\maketitle

\begin{abstract}

{\swr
{\em Zombie fires} in peatlands disappear from the surface, smoulder underground during the winter, and `come back to life' in the spring. They can release hundreds of megatonnes of carbon into the atmosphere per year and are believed to be caused by surface wildfires.
}
Here, we propose rate-induced tipping (R-tipping) to a {\em  subsurface hot metastable state} in bioactive peat soils as a main cause of Zombie fires. Our hypothesis is based on a conceptual soil-carbon model
subjected to realistic changes in weather and climate patterns, including global warming scenarios and 
summer heatwaves.

Mathematically speaking, R-tipping to the hot metastable state is a nonautonomous  instability, due to crossing an elusive {\em quasithreshold}, in a multiple-timescale dynamical system. To explain this instability, we provide a 
framework combining a special compactification technique with concepts from geometric singular perturbation theory. 
This framework allows us to reduce an R-tipping problem due to crossing a quasithreshold to a heteroclinic orbit problem in a singular limit. We
 identify generic cases of 
 {\swb tracking-tipping transitions}
 via: (i) unfolding of a  codimension-two {\em heteroclinic folded saddle-node  type-I singularity} for global warming, and (ii) analysis of a codimension-one {\em saddle-to-saddle hetroclinic orbit} for summer heatwaves, in turn revealing new types of excitability quasithresholds.
 
\end{abstract}


\section{Introduction}

Tipping points, or critical transitions, are instabilities known to occur in natural and human systems subjected to changing external conditions, or external inputs. They may be explained in layman's terms as 
{\sw an abrupt}
and large change in the state of a system in response to a small or slow change in the external input. The change in the state of the system may be permanent or transient.

{\sw Climate change is an important factor for tipping points in natural systems. }
An immense amount of research is being conducted to predict and prevent its worst effects. Hand in hand with this research effort, interest in tipping has accelerated due to theorised~\cite{Lenton1786,ritchie_overshooting_2021} and observed~\cite{caesar_rahmstorf_robinson_feulner_saba_2018, boers_rypdal2021} tipping points in  the  earth system caused by climate change.
Many of the tipping elements  identified in~\cite{Lenton1786}, for example the loss of  Arctic sea ice  or the shutdown of the Atlantic Meridional Overturning Circulation, can be captured by elaborate, high-resolution mathematical models referred to as General Circulation Models (GCMs). 
In this paper, we use a conceptual soil-carbon model with time-varying climate as an external input to describe 
a tipping element that is not captured by GCMs:
A release of gigatonnes of carbon from temperature-sensitive peat soils into the atmosphere via so-called ``Zombie fires"~\cite{BBC2021} that disappear from the surface, smoulder underground during the winter, and ``come back to life" in the spring.

Owing to the explicit time dependence of the external input, the ensuing dynamical system is nonautonomous.
This means that analysis of tipping points requires, in general, techniques beyond classical autonomous stability theory~\cite{ashwin2012,okeeffe_tipping_2020,Wieczorek2021R}.
Nonetheless, it is useful to consider the corresponding {\em autonomous frozen system} with fixed-in-time inputs. In the frozen system we identify a desired stable state, and refer to this state as the {\em base state}.
When the external input changes over time,  the shape and position of the base state may change too, and the nonautonomous system will try to {\em track} the moving base state. However, 
sometimes tracking is not possible and tipping occurs.
For example, the base state may lose stability or disappear in a classical bifurcation at some {\em critical level} of the input. If this bifurcation is {\em dangerous}~\cite{thompson_sieber_2011}, the system tips to a different stable state, referred to as an {\em alternative stable state}.
We then say the system undergoes {\em bifurcation-induced tipping}, or in short {\em B-tipping}~\cite{thompson_sieber_2011,kuehn2011,ashwin2012}. 
Another, arguably more interesting example is when the external input changes faster than some {\em critical rate}, the nonautonomous system deviates too far from the changing base state, crosses some {\em threshold}
or {\em quasithreshold}~\cite{fitzhugh_1955}
and tips to an {\em alternative state}. 
Such tipping is caused entirely by the {\em rate of change} of the external input and 
we say the system undergoes {\em rate-induced tipping}, or in short {\em R-tipping}~\cite{ashwin2012,Wieczorek2021R}.
%
Crucially, unlike B-tipping, R-tipping can occur to an {\em  alternative transient state} that lasts for a finite time\footnote{This is in contrast to an {\em alternative stable state} that lasts forever.}, after which the system returns to the {\sw base state~\cite{wieczorek_excitability_2011,Mitry2013,Vanselow_2019,hastings2018}}. Such R-tipping is referred to as ``reversible" in~\cite{Wieczorek2021R}. Systems that exhibit reversible R-tipping are also known as {\em excitable systems}~\cite{Izhikevich2006}.
%
The peat soil instability studied in this paper is an example of
{\em reversible R-tipping} to an {\em alternative metastable state} that lasts for a long but finite time, that occurs due to crossing an elusive quasithreshold.

\begin{figure}[t]
    \centering
    \includegraphics[width=0.8\textwidth]{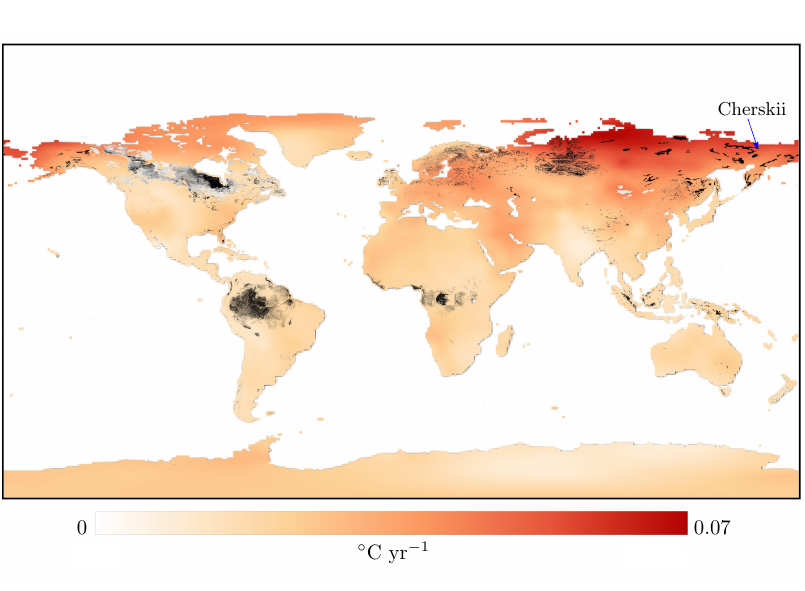}
    \caption{
     (Colourscale)  Historical global rates of warming from the period 1951-1980 to the period 2017-2021  obtained using  observations compiled by Berkeley Earth~\cite{berkely_earth_data,Berkely_Earth_Temp_generation}, together with (black and gray) the global distribution of peat soils obtained from PEATMAP~\cite{PEATMAP2018}. Note the largest warming rates at higher latitudes, in areas with a significant concentration of peat soils, such as (blue arrow) Cherskii in Northern Siberia.
   }
    \label{fig:RateMap}
\end{figure}

The main motivation for our study is a combination of two environmental  features of the Arctic. First is the organic carbon content.
Estimates for the organic carbon contained in Northern and Arctic permafrost peat soil alone range from approximately 500Gt~\cite{bg-9-4071-2012} to approximately 1700Gt~\cite{tarnocai_soil_2009}. 
For comparison, the atmospheric carbon pool is estimated to be approximately 850 Gt~\cite{ClimateChange_ThePhysicalScienceBasis}.
Second is the rate of atmospheric warming and the increasing trends in summer heatwaves.
Due to so called ``Arctic Amplification''~\cite{screen_simmonds_2010},  both the so-far observed and future predicted warming for Arctic regions are approximately double the global mean~\cite{overland_urgency_2019}. 
This means that the Arctic is home to massive deposits of ancient peat carbon and is the fastest warming region on the planet.
To visualise this combination of environmental features we combine in fig.~\ref{fig:RateMap} (colourscale) recent rates of global warming\footnote{Note that these recent short-term rates of global warming exceed the long-term rates in the CMIP5 outputs.} 
and (greyscale) the global distribution of peat soils.
%
In addition to the increase in the mean global temperature, there is an increase in the intensity, frequency and duration of summer heatwaves in the Northern Hemisphere~\cite{perkins-kirkpatrick_lewis_2020}, 
with the Arctic temperature record a scorching $38^\circ$C in 2020~\cite{united_nations_2021}.
Since peat soils are bioactive and thus temperature sensitive~\cite{Makiranta2009}, such conditions
can lead to 
{\sw {\em thermal runaway} in the soil}. This is the reason why northern latitude peat-soils were identified as a potential tipping element in~\cite{Lenton1786}.
To the best of our knowledge, the first example of 
R-tipping in peat soils: 
{\sw thermal runaway}
triggered by {\sw the rate of} 
atmospheric warming, 
was reported, but not emphasised, by Khvorostyanov et al.~\cite[Fig.4(a)]{khvorostyanov_vulnerability_2008-1}. 
Later, Luke and Cox~\cite{luke_soil_2011} proposed a conceptual soil-carbon model that exhibits a short-lived explosive release of soil carbon into the atmosphere above some critical rate of global warming, which they termed the ``compost-bomb instability"; see also~\cite{Clarke2021}. The dynamical mechanism responsible for this {\sw R-tipping} instability was explained by Wieczorek et al.~\cite{wieczorek_excitability_2011}.
Additionally, it has been known that spontaneous combustion is the main cause of fires at composting facilities~\cite{Buggeln2002}, which can then spread, e.g. the recent Wennington fire in London~\cite{sawer_2022}.
%

In the first, mathematical modelling part of the paper in Sections~\ref{sec: rtscm} and~\ref{sec:RHMSZF}, we:
\begin{itemize}
    \item 
    Modify the  conceptual 
    model introduced by Luke and Cox~\cite{luke_soil_2011} with a more realistic microbial soil respiration function.
    \item
    Show that the modified  conceptual  model has a new {\em alternative hot metastable state} and reproduces the key features of the medium-complexity model from~\cite[Fig.4(a)]{khvorostyanov_vulnerability_2008-1}. 
    \item
    Demonstrate R-tipping to the hot metastable state in the modified conceptual model for realistic climate change scenarios including global warming and summer heatwaves.
    \item
    Based on the above, propose an explanation for ``Zombie fires" in peatlands~\cite{BBC2021,Scholten2021,mccarty_arctic_2020,vaughan_arctic_2020}.
\end{itemize}
In the second, mathematical analysis part of the paper in Sections~\ref{sec: mtc}--\ref{sec:HSA}, we identify non-obvious dynamical mechanisms that are responsible for the R-tipping instability to the hot metastable state. The main obstacles to
analysis of this instability are twofold: the conceptual soil-carbon model is a nonautonomous dynamical system, so it does not have any equilibria or compact invariant sets, and the R-tipping is due to crossing an elusive quasithreshold in the phase space~\cite{Mitry2013,perryman_adapting_2014}.
To overcome these obstacles we combine three different strategies.
Firstly, we consider external inputs that decay to a constant at infinity. Secondly, we compactify the problem to include the equilibrium base states 
for the autonomous limit systems from infinity~\cite{Wieczorek_2021}. 
These two strategies alone work for R-tipping due to crossing  regular thresholds that
are anchored at infinity by  unstable compact invariant sets called {\em regular R-tipping edge states}~\cite[Sec.4]{Wieczorek2021R}. However,  quasithresholds do not contain such edge states~\cite[Sec.8]{Wieczorek2021R}. Therefore, thirdly, 
we exploit
large timescale separation in the 
model and apply concepts from geometric singular perturbation theory {\sw (GSPT)}. 
Specifically, we:
\begin{itemize}
    \item
    Define R-tipping due to crossing a quasithreshold in the nonautonomous system in terms of {\em slow manifolds} and {\em canard trajectories}~\cite{Wechselberger2013} for the autonomous compactified system. 
    \item
    Identify {\em singular R-tipping edge states}: special points called {\em folded singularities}~\cite{SZMOLYAN2001419} that arise in the reduced (slow) system, and {\em new saddle equilibria} that arise in the layer (fast) system. 
    \item
    Reduce an R-tipping problem due to crossing a quasithreshold to a heteroclinic orbit connecting the base state for the past limit system to a singular R-tipping edge state. We then use this reduction to identify four different cases of 
    {\swb tracking-tipping transition:}
    \begin{itemize}
        \item
        For global warming, three (slow) cases
        are identified via unfolding of a  codimension-two {\em  non-central heteroclinic folded saddle-node type-I singularity} that arises in the reduced system.
        \item
        For a summer heatwave,  a fourth (fast) case 
        is identified via analysis of a  codimension-one {\em heteroclinic orbit} connecting the base state for the past limit system to a new saddle  equilibrium that arises in the layer system.
        \end{itemize}
    \item
    Show that a quasithreshold gives rise to {\em critical ranges of rate of change} of the external input rather than isolated critical rates. 
    Furthermore, we reveal new types of quasithresholds  through analysis of canard trajectories associated with singular R-tipping edge states.
\end{itemize}

\section{The Nonautonomous Soil-Carbon Model}
\label{sec: rtscm}

The starting point of our analysis is a 
discussion of the conceptual soil-carbon model introduced
by Luke and Cox~\cite{luke_soil_2011}.
This model describes the time evolution of the soil temperature $T$ and soil carbon concentration $C$ in peat soils
\begin{align}
\label{eq: Intro Model}
    \mu\frac{dT}{dt} &= -\lambda(T-T_{a}(rt)) + A\, C\,R_s(T),\\
\label{eq: Intro Model2}
    \frac{dC}{dt} &= \Pi-C\,R_s(T).
\end{align}
The model incorporates three soil processes and one time-varying external input.
The first process describes {\sw balancing of the soil, $T$, and atmospheric, $T_a$, temperatures towards a thermal equilibrium}
according to Newton's Law of Cooling,
at a rate that depends on the soil-to-atmosphere heat transfer coefficient $\lambda$ and the specific heat capacity of the soil $\mu$.
The second process describes a linear increase in the soil carbon concentration $C$ over time at a
rate $\Pi$, due to carbon generated from decaying plant litter and other processes referred to as  ``Gross Primary Production"~\cite{roy_terrestrial_2001}.\footnote{Note that $dC/dt$ is the ``Net Primary Production" in the model.} 
The third and only nonlinear process in the model describes temperature-sensitive microbial activity in the soil in terms of the soil respiration function $R_s(T)$. This
process couples the dynamics of $T$ and $C$, and is discussed in detail in section~\ref{sec: rtscm}\ref{sec:msrhms}. 
Atmospheric temperature $T_a(rt)$ is a time-varying external input, which represents weather anomalies or climate variation.
The {\em rate parameter} $r > 0$ quantifies the time scale of climatic variability and is the key {\em input parameter} in the model.
An important aspect of our study is that we use realistic values of the soil parameters based on~\cite{luke_soil_2011}, which are given in Table 1
in the  electronic appendix,
a realistic soil respiration function $R_s(T)$ introduced in 
section~\ref{sec: rtscm}\ref{sec:msrhms}, and realistic climate-change scenarios $T_a(rt)$ based on real weather data from Cherskii in Siberia~\cite{noauthor_chersky_nodate}; see the arrow in fig.~\ref{fig:RateMap}.

To facilitate {\sw the} analysis, we introduce  a {\em small parameter} $\epsilon = \mu/A \ll 1$ 
and rewrite {the soil-carbon model~\eqref{eq: Intro Model}--\eqref{eq: Intro Model2} as a fast-slow nonautonomous dynamical system}~\cite{wieczorek_excitability_2011}:
\begin{align} \label{eq: SmallParam}
    \epsilon\frac{dT}{dt} & = f_1(T,C,T_a(rt)) := -\frac{\lambda}{A}(T-T_{a}(rt)) + C\,R_s(T),\\
    \label{eq: SmallParam2}
    \frac{dC}{dt} &=  f_2(T,C)\quad\quad\quad\, := \Pi-C\,R_s(T),
\end{align}
where we define $f_1$ and $f_2$ for convenience.
Note that system~\eqref{eq: Intro Model}--\eqref{eq: Intro Model2} or~\eqref{eq: SmallParam}--\eqref{eq: SmallParam2} does not have any equilibria (stationary solutions) owing to the time-varying external input $T_a(rt)$.

\begin{figure}[t]
    \centering
    \includegraphics[width=0.45\textwidth]{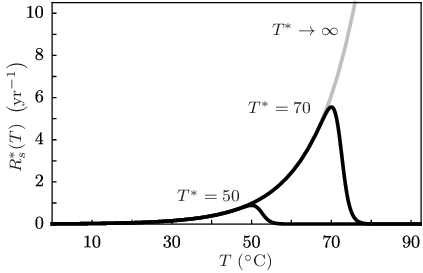}
    \caption{(Black) The modified non-monotone  microbial soil respiration function $R_s^{*}(T)$
    in~\eqref{eq:new_r} with $c=10$ and the die-off temperature $T^{*}=50^\circ$C and $70^\circ$C. (Grey)
    In the limit $T^{*}\to\infty$, we recover the unmodified monotone  microbial soil respiration function
    $R_s^{\dagger}(T)$ in~\eqref{eq:old_r}. See Table 1 in the electronic appendix for other parameter values.
    }
    \label{fig:Rs Comparison}
\end{figure}

\subsection{ Modified Microbial Respiration}
\label{sec:msrhms}

Our contribution to the model introduced
by Luke and Cox is a modification of the 
microbial soil respiration function from~\cite{luke_soil_2011}.
At low to moderate soil temperatures $T$, microbial soil respiration can be described by a
$Q_{10}$ exponential function~\cite{kirschbaum_temperature_1995}: 
\begin{equation}
\label{eq:old_r}
R_s^{\dagger}(T)= R_{s}(0)\,e^{\alpha T} \quad \mbox{with} \quad  \alpha = \ln{(Q_{10})}/10,
\end{equation}
where the dimensionless, soil-specific parameter $Q_{10}$ may be estimated from experimental data~\cite{kirschbaum_temperature_1995}.
While the monotone $R_s^{\dagger}(T)$ in~\eqref{eq:old_r} captures 
{\sw thermal runaway -} the  
{\sw process responsible}
 for the R-tipping instability, 
it becomes unrealistic at high $T$. 
Specifically, $T$ quickly increases to unrealistically high levels due to microbial soil respiration alone~\cite{luke_soil_2011,wieczorek_excitability_2011}.
To address this issue, 
we account for an important limitation, that is, soil microbes die above some {\em die-off temperature} $T=T^*$. Specifically, we construct a modified non-monotone microbial soil respiration function 
\begin{equation}
\label{eq:new_r}
R^*_{s}(T) = {R}_{s}(0)\,\frac{e^{\alpha b} + e^{-c\alpha b}}{e^{-\alpha (T-b)} + e^{c\alpha (T-b)}}
\quad \mbox{with} \quad
b = T^{*}+\frac{\ln{c}}{\alpha+c\alpha},
\end{equation}
that agrees with the $Q_{10}$ exponential growth~\eqref{eq:old_r} for $T < T^*$,  has a maximum at $T = b \approx T^*$, and decays exponentially for $T > T^*$; see fig.~\ref{fig:Rs Comparison}.
Such $R^*_{s}(T)$ {\sw  captures thermal runaway and, additionally, stops it at high but realistic levels of  $T$}. 
{\sw This} construction introduces 
two additional parameters, namely the die-off temperature $T^*$, and the ratio $c$ of the exponential decay rate for $T>T^*$ and exponential growth rate $\alpha$ for $T<T^*$.
In the remainder of the paper, we use
$$
R_s(T) = R_s^*(T).
$$

\subsection{The  Autonomous Frozen System and the Moving Equilibrium}
\label{sec:qse}

To gain insight into the tipping mechanisms in the nonautonomous system~\eqref{eq: SmallParam}--\eqref{eq: SmallParam2}, or equivalently system~\eqref{eq: Intro Model}--\eqref{eq: Intro Model2},
we set $r=0$ and consider properties of the resulting autonomous {\em frozen system} 
with fixed-in-time $T_a$.
The frozen system has just one equilibrium
\begin{align}
\label{eq:qse}
e(T_a) &=  (T^{e}(T_a),C^{e}(T_a)) =
\left(T_a+\frac{A\Pi}{\lambda},\frac{\Pi}{R_s(T_a+\frac{A\Pi}{\lambda})}\right),
\end{align}
which is the ``base state" described in the introduction.
The position of $e(T_a)$  in the  phase plane $(T,C)$ changes with $T_a$, but the equilibrium
remains  linearly stable 
and globally attractive (attracts all initial conditions)
within the realistic 
range of  $T_a$ used in our study.
Thus, we can exclude the possibility of B-tipping from $e(T_a)$~\cite{ashwin2012}.
 Next, we consider the stable equilibrium of the frozen system parameterised by time $t$ for a given input $T_a=T_a(rt)$, which we 
denote 
\begin{equation}
    \label{eq:me}
e(T_a(rt)),
\end{equation}
and refer to as the
{\em moving  stable equilibrium}~\cite{okeeffe_tipping_2020,Wieczorek2021R}.
The moving stable equilibrium is not a solution to the nonautonomous
system~\eqref{eq: SmallParam}--\eqref{eq: SmallParam2}. However, it can serve as a useful point of reference.
We follow the approach of~\cite{ashwin2017,okeeffe_tipping_2020,Wieczorek2021R} and relate solutions of the nonautonomous system~\eqref{eq: SmallParam}--\eqref{eq: SmallParam2} to 
$e(T_a(rt))$ for different rates $r>0$.
For small $r$, solutions to~\eqref{eq: SmallParam}--\eqref{eq: SmallParam2} started near $e(T_a(rt))$
are guaranteed to stay near or track $e(T_a(rt))$~\cite[Th.7.1]{Wieczorek2021R}. However, a nonautonomous R-tipping instability in the form of a  large  transient departure from $e(T_a(rt))$ 
may appear for larger $r$. 
\begin{figure}[t]
    \centering
    \includegraphics[width=0.8\textwidth]{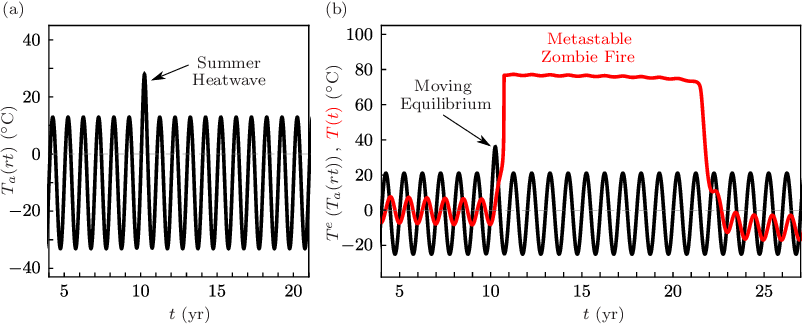}
    \caption{
    (a)
      Seasonal variations of atmospheric temperature $T_a(rt)$ with a summer heatwave
      in year 10, based on observations at Cherskii in Northern Siberia \cite{noauthor_chersky_nodate}.
     (b) Time evolution of (black) moving equilibrium soil temperature $T^{e}(T_a(rt))$ from~\eqref{eq:qse} and (red) the actual soil temperature $T(t)$ obtained by solving~\eqref{eq: SmallParam}--\eqref{eq: SmallParam2} with
     the input $T_a(rt)$ from panel (a). The model is initialised at $t=-10\,$yr, and
    initial soil temperature is $T(-10)= -10^\circ$C. 
     Initial soil carbon $C(-10) = 120$\,kg\,m$^{-2}$ corresponds to $\approx3.64$\,m soil depth
    assuming a volumetric carbon density of $33$\,kg\,m$^{-3}$~\cite{khvorostyanov_vulnerability_2008-1, zimov_permafrost_2006}.
    $R_s(T)=R_s^*(T)$ given in~\eqref{eq:new_r} with $T^* = 70^\circ$C.
    See the electronic appendix for  other parameter values in Table 1 and  details of $T_a(rt)$ in section 1.3.
    }
    \label{fig: Periodic With Hot Summer}
\end{figure}
\section{ R-tipping to a Hot Metastable State and Zombie Fires}
\label{sec:RHMSZF}

The more realistic $R_s(T)=R_s^*(T)$ introduced in section~\ref{sec: rtscm}
gives rise to a {hot metastable state} at high soil temperatures $T \approx T^*$. This state is reminiscent of ``Zombie fires" observed in tropical and arctic peatlands~\cite{BBC2021,Scholten2021,mccarty_arctic_2020,vaughan_arctic_2020}:
such fires appear to be extinguished, but smoulder underground throughout the winter and re-emerge the following year~\cite{Scholten2021}.
Zombie fires are generally believed to happen as a result of surface wildfires~\cite{BBC2021}. However, as far as we know, there have been no attempts to explain this phenomenon. Here,
we propose a hypothesis that R-tipping to a subsurface hot metastable state due to atmospheric warming is a main cause of Zombie fires. Our hypothesis is based on two remarkable results of the soil-carbon model~\eqref{eq: SmallParam}--\eqref{eq: SmallParam2} for realistic soil parameters
and different climate-change scenarios $T_a(rt)$ based on real weather data~\cite{noauthor_chersky_nodate}.

\begin{figure}[t]
    \centering
    \includegraphics[width=0.8\textwidth]{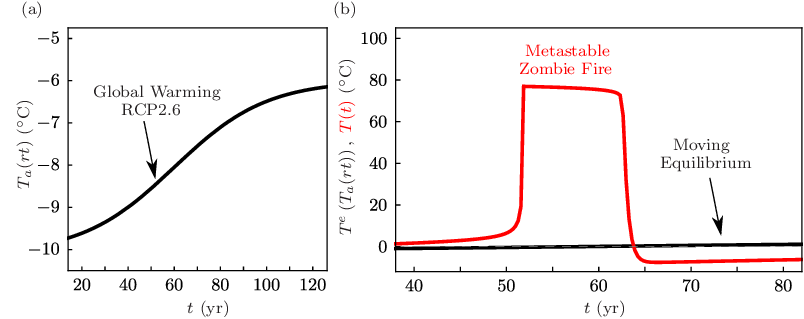}
    \caption{
    (a)
    A realisation 
    within the range of uncertainty of
    the ``very stringent" low-emissions RCP2.6 global warming scenario
    for Arctic regions~\cite{pachauri_climate_2015, overland_urgency_2019},
    based on observations at Cherskii in Northern Siberia~\cite{noauthor_chersky_nodate}.
    (b)
    Time evolution of (black) the moving equilibrium soil temperature $T^{e}(T_a(rt))$ from~\eqref{eq:qse} and (red) the actual soil temperature $T(t)$ obtained by solving~\eqref{eq: SmallParam}--\eqref{eq: SmallParam2} 
    with the input $T_a(rt)$ from panel (a). 
    Initial soil temperature is $T(0)= T_a^-+\frac{A\Pi}{\lambda}\approx-1.95^\circ$C. 
    Initial soil carbon $C(0)=120$\,kg\,m$^{-2}$ and $R_s(T)$ is 
     the same as in fig.~\ref{fig: Periodic With Hot Summer}.
     See the electronic appendix for  other parameter values in Table 1 and  details of $T_a(rt)$ in section 1.4.
    }
    \label{fig: Climate Change}
\end{figure}

{\bf R-tipping to a hot metastable state in the 
conceptual
model~\eqref{eq: SmallParam}--\eqref{eq: SmallParam2}
can be triggered by realistic climate patterns ranging from summer heatwaves to global warming scenarios.} 
Figure~\ref{fig: Periodic With Hot Summer}
shows the soil temperature change in response to seasonal variations of the atmospheric temperature. Rather amazingly, {\em a summer heatwave} in year 10 breaks the seasonal response pattern and  triggers a sudden transition to a hot metastable state
that lasts for over a decade and releases 
most of the carbon from the soil into the atmosphere. Following this, the system settles to a lower than initial soil temperature pattern for a refractory period of over a century, during which time both soil carbon and soil temperature slowly increase back to their initial seasonal patterns.
Figure~\ref{fig: Climate Change} shows the mean soil temperature change in response to  a slow increase in the mean atmospheric temperature of $4^\circ$C over 200 years.\footnote{We use $4^\circ$C in this example because the global climate change mitigation target is an increase in the global mean temperature of $2 ^\circ$C by 2100 compared to pre-industrial levels (1850–1900)~\cite{IPCCsr15, unfccc2015adoption}. 
The $2 ^\circ$C target is within the range of uncertainty of CMIP5 outputs under two greenhouse gas concentration pathways: the ``very stringent'' RCP2.6 that gives a $1.5^\circ$C increase, and the ``intermediate" RCP4.5 that gives a $2.4^\circ$C increase~\cite{overland_urgency_2019}. Due to the well documented phenomenon of ``Arctic amplification''~\cite{overland_urgency_2019,screen_simmonds_2010}, this target corresponds to an increase of roughly $4^\circ$C in the Arctic mean temperature over the same period.}
In this scenario, a sudden transition to a hot metastable state that lasts for over a decade is triggered, and astonishingly occurs after only a modest mean temperature increase of $\approx1.1^\circ$C.
Following this, the system settles to a lower than initial soil temperature for a refractory period of over a century, during which time both soil carbon and soil temperature slowly increase to their moving equilibrium levels.

\begin{figure}[t]
   \centering
    \includegraphics[width=0.7\textwidth]{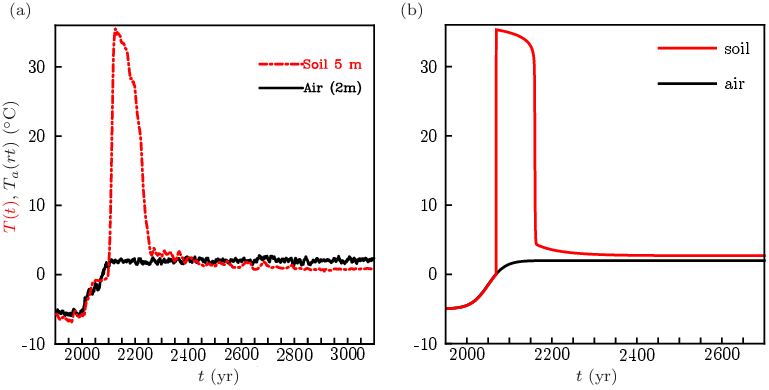}
    \caption{(a) The output of a PDE model~\cite{khvorostyanov_vulnerability_2008} parameterised by data from near Cherskii in Siberia to
    simulate a soil column with a depth of 12m and a volumetric carbon density of $33$\,kg\,m$^{-3}$, reproduced from~\cite[Fig.4(a)]{khvorostyanov_vulnerability_2008-1}. 
    Shown are time evolution of (black) air temperature at 2 metres above the ground, and (red) soil temperature at a depth of 5m.
    (b)
    Time evolution of (red) the soil temperature {$T(t)$} obtained by solving~\eqref{eq: SmallParam}--\eqref{eq: SmallParam2} 
    with the (black) input $T_a(rt)$.
    For comparison with (a),
    we use $\Pi=0.09$ kg m$^{-2}$ yr$^{-1}$,  $R_s(T)=0$ for $T\leq0$ as in~\cite{khvorostyanov_vulnerability_2008-1}, and $R_s(T)=R_s^*(T)$  given in~\eqref{eq:new_r} with $T^* = 30^\circ$C for $T>0$.
    Initial soil temperature is $T(0) = -6^\circ$C. Initial soil carbon is $C(0) = 396$\,kg\,m$^{-2}$, 
    which corresponds to a soil column with a depth of $12$\,m and a volumetric carbon density of $33$\,kg\,m$^{-3}$ matching~\cite{khvorostyanov_vulnerability_2008-1}.
    See the electronic appendix for  other parameter values in Table 1 and  details of $T_a(rt)$ in section 1.5.
    }
    \label{fig: Khvorostyanov Comparison}
\end{figure}

{\bf R-tipping to a hot metastable state  in the conceptual
model~\eqref{eq: SmallParam}--\eqref{eq: SmallParam2} shows qualitative and quantitative agreement with the results of intermediate-complexity PDE models.} 
The conceptual Ordinary Differential Equation (ODE) model~\eqref{eq: SmallParam}--\eqref{eq: SmallParam2} with the more realistic $R_s(T)=R_s^*(T)$
captures the key nonlinearities of the soil-carbon system. This is evidenced by the ODE model
reproducing both qualitatively and quantitatively the peat soil
instability predicted by a medium-complexity 
Partial Differential Equation (PDE) model of Siberian permafrost carbon dynamics under climate change~\cite{khvorostyanov_vulnerability_2008-1,khvorostyanov_vulnerability_2008}.
Specifically, fig.~\ref{fig: Khvorostyanov Comparison} 
shows that the ODE model manages to capture the four key features of the PDE dynamics: permafrost thawing, R-tipping to the hot metastable state that lasts for half a century, followed by a sudden soil cooling to slightly above the air temperature. 
%

\section{The Multiscale Autonomous Compactified System
}
\label{sec: mtc}

To analyse and understand
non-obvious dynamical mechanisms that are responsible for the R-tipping instabilities in figs.~\ref{fig: Periodic With Hot Summer} and~\ref{fig: Climate Change}, we:
\begin{itemize}
    \item[$\bullet$]
    Consider 
    external inputs $T_a(rt)$ that 
    tend exponentially to a constant at infinity.
    \item[$\bullet$]
    Reformulate the ensuing two-dimensional nonautonomous  system~\eqref{eq: SmallParam}--\eqref{eq: SmallParam2} as a three-dimensional autonomous compactified system~\cite{Wieczorek_2021}.
    \item[$\bullet$]
    Identify three different timescales in the soil-carbon system.
    
\end{itemize}


We choose to work with  external inputs $T_a(rt)$ 
that  decay exponentially to a constant $T_a^\pm$ as time $t$  tends to $\pm\infty$. To be precise, we follow~\cite[Def.6.1]{Wieczorek2021R} and 
\begin{defn}
\label{defn:expbiconst}
We say $T_a(rt)$ is {\em exponentially bi-asymptotically constant} if,
for all $r > 0$,
\begin{align*}
    \lim_{t\to\pm\infty}T_a(rt)=T_a^\pm\in\mathbb{R}\;\;\;
    \mbox{and}\;\;
    \lim_{t\to\pm\infty}\frac{dT_a(rt)/dt}
{e^{\mp\rho rt}}\in\mathbb{R},
\end{align*}
 for some {\em decay coefficient} $\rho>0$.
\end{defn}
Thus, we can define the autonomous {\em past limit system} with $T_a(rt)=T_a^-$, 
\begin{align}
\label{eq:ls-}
    \epsilon\, \frac{dT}{dt} = f_1(T,C,T_a^-),\quad\quad
    \frac{dC}{dt} = f_2(T,C),
\end{align}
and the autonomous {\em future limit system} with $T_a(rt)=T_a^+$,
\begin{align}
    \label{eq:ls+}
    \epsilon\, \frac{dT}{dt} = f_1(T,C,T_a^+),\quad\quad
    \frac{dC}{dt} = f_2(T,C), 
\end{align}
which are examples of the frozen system. 
We note that, unlike the nonautonomous system~\eqref{eq: SmallParam}--\eqref{eq: SmallParam2}, the autonomous past~\eqref{eq:ls-} and future~\eqref{eq:ls+} limit systems contain the equilibria 
\begin{equation}
\label{eq:limeq}
e^- := e(T_a^-)\quad\mbox{and}\quad e^+ :=e(T_a^+),
\end{equation}
respectively, and
$e(T_a(rt))\to e^\pm$ as $t\to\pm\infty$ for any $r>0$.

In the next step, we include these limit systems and their equilibria in the model; 
{\sw see}~\cite{Wieczorek2021R,Wieczorek_2021} and section~2 of the electronic appendix for full details.
We introduce a bounded dependent variable
\begin{align}
\label{eq:comptrans}
    s = g_{\nu}(rt) = \tanh\left(\frac{\nu}{2}\,rt\right)\in(-1,1),
\end{align}
and reformulate the two-dimensional nonautonomous system~\eqref{eq: SmallParam}--\eqref{eq: SmallParam2} on $\mathbb{R}^2$
as a three-dimensional {\em autonomous compactified system} on the {\em extended phase space} $\mathbb{R}^2\times[-1,1]$,
\begin{align}
\label{eq: CompactEqns}
        \epsilon\; \frac{dT}{dt} &= f_1(T,C,T_a^{\nu}(s)),\\
    \label{eq: CompactEqns2}
        \frac{dC}{dt} &= f_2(T,C),\\
     \label{eq: CompactEqns3}
        \frac{1}{r}\;  \frac{ds}{dt} &= \frac{\nu}{2}(1-s^2),
\end{align}
with the continuously extended external input
\begin{equation}
\label{eq: LambdaGeneral}
    T_a^{\nu}(s) = 
    \left\{
    \begin{array}{rcl}
    T_a\left( g^{-1}_{\nu}(s) \right) & \text{for} & s\in(-1,1),\\
    T_a^+ & \text{for} & s=1,\\
    T_a^- & \text{for} & s=-1,
    \end{array}
  \right.
\end{equation}
and {\em compactification parameter} $\nu$.
A particular advantage of compactification is that the flow-invariant planes 
\begin{align}
    \label{eq:inv_planes}
\mathbb{R}^2\times\{-1\}
\quad\mbox{and}\quad
\mathbb{R}^2\times\{1\},
\end{align}
of the extended phase space contain equilibria $e^{-}$ and $e^{+}$ of the autonomous past~\eqref{eq:ls-} and future~\eqref{eq:ls+} limit systems, respectively.
When embedded in the extended phase space, 
$e^-$ gains one unstable eigendirection with positive eigenvalue $\nu\, r > 0$ and becomes a hyperbolic saddle,
\begin{align}
\label{eq:limeqcomp-}
&\tilde{e}^- = {\left(e^-,-1 \right)},
\end{align}
 whereas $e^+$ gains one additional stable eigendirection with negative eigenvalue $-\nu\, r < 0$
and becomes a hyperbolic sink, 
{\sw which is the only attractor for the compactified system~\eqref{eq: CompactEqns}--\eqref{eq: CompactEqns3},}
\begin{align}
\label{eq:limeqcomp+}
\tilde{e}^+ = {\left(e^+,+1 \right)}.
\end{align}
%
Furthermore, we note that the moving equilibrium $e(T_a(rt))$ with $t\in\mathbb{R}$ corresponds to
$$
\tilde{e}(s) := \left(
e(T_a\left(g^{-1}_{\nu}(s)\right),s
\right)~\mbox{with}~s\in(-1,1),
$$
in~\eqref{eq: CompactEqns}--\eqref{eq: CompactEqns3}, and $\tilde{e}(s)\to\tilde{e}^\pm$ as $s\to\pm 1$.

The left-hand side of the compactified system~\eqref{eq: CompactEqns}--\eqref{eq: CompactEqns3} shows that the soil-carbon system
may evolve on up to three different timescales, depending on the rate parameter $r$.
The {\em slow time} $t$ is the timescale of the {\em slow variable} $C$. The {\em fast time} $\tau=t/\epsilon$
is the timescale of the {\em fast variable} $T$. The {\em third time} $u = rt$ is the timescale of the external input $T_a(rt)$, which is represented by the additional variable $s$.
The magnitude of the third time $u$ relative to $t$ and $\tau$ depends 
on the magnitude of the rate parameter $r$.
Here, we consider different two-timescale limits of this three-timescale problem. 
{\sw
A more unifying view of the system dynamics can be obtained through analysis of the three-timescale problem~\cite{cardin2017}~[Sec.6, electronic appendix], which is beyond the scope of this paper.
}

\section{Defining R-tipping due to Crossing a Quasithreshold}
\label{sec:frozen}

The multiple-timescale soil-carbon system
can be viewed as a singular perturbation problem~\cite{kuehn_2016}. 
In this section, we combine compactification with concepts and techniques of {\sw GSPT} to:
\begin{itemize}
    \item 
    Define R-tipping due to crossing a quasithreshold in the 
    soil-carbon
    system~\eqref{eq: SmallParam}--\eqref{eq: SmallParam2}.
    \item
    Give intuition for when to expect such R-tipping. 
\end{itemize}

\begin{figure}[]
    \centering
   \includegraphics[width=0.8\textwidth]{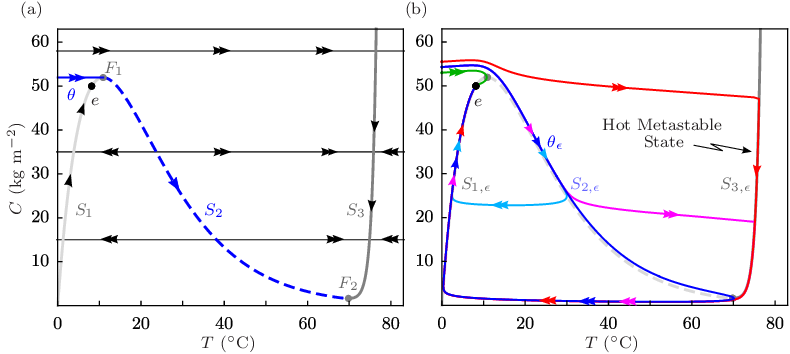}
   \caption{(a) The critical manifold $S(T_a)$ with (grey) two normally hyperbolic attracting branches, $S_1(T_a)$ and $S_3(T_a)$,  and (dashed blue) a normally hyperbolic repelling branch $S_2(T_a)$, for the frozen system with $T_a = 0$ and $\epsilon = 0$. 
   Shown also are {(single arrows) 
   {\sw slow timescale trajectories}
   of the reduced problem~\eqref{eq: frozenreduced} along $S(T_a)$, 
   and  (double arrows) 
   {\sw fast timescale trajectories}
   of the layer problem~\eqref{eq: frozenLayer}
   for different values of $C$.} 
   (b) Five {\sw fast-slow composite} trajectories in 
   the autonomous frozen system with $T_a = 0$, $\epsilon \approx 0.064$, the same initial soil temperature $T(0)=0$, and different initial carbon levels: (green) $C(0) = 53$, (cyan) $C(0) = 54.287996178057110079$, (blue) $C(0) \approx 54.28799617805711008$, (magenta) $C(0) = 54.287996178057110088$, and
   (red) $C(0) = 55.5$. 
   Included for reference is (grey) $S(T_a)$.
   See Table 1 in the electronic appendix for other parameter values.
   }
    \label{fig:frozensystem}
\end{figure}

It is convenient to start the discussion with the autonomous frozen system  obtained by setting $r=0$ in system~\eqref{eq: SmallParam}--\eqref{eq: SmallParam2}, so that $T_a$ becomes a fixed-in-time input parameter. The frozen system is a 1-fast 1-slow singular perturbation problem with a small parameter $\epsilon$. Taking the limit $\epsilon\to0$ in the slow time $t$ gives 
a singular {\em reduced problem}
\begin{equation}
\begin{split}
\label{eq: frozenreduced}
    \frac{dC}{dt} &= f_2(T,C),
\end{split}
\end{equation}
{\sw with {\em slow timescale solutions}} restricted to the  one-dimensional critical manifold
\begin{align*}
    S(T_a) = \left\{(T,C)~:~f_1(T,C,T_a) = 0 \right\}\subset\mathbb{R}^2.
\end{align*}
Taking the limit $\epsilon\to0$ in the fast time $\tau=t/\epsilon$ gives 
a {\sw regular} {\em layer problem}
\begin{align}
\begin{split}
\label{eq: frozenLayer}
    \frac{dT}{d\tau} &= f_1(T,C,T_a),
\end{split}
\end{align} 
{\sw
with {\em fast timescale solutions} for a fixed-in-time $C$.
}
%
Note that, for a given $T_a$, the critical manifold $S(T_a)$ consists of all branches of equilibria of the layer problem~\eqref{eq: frozenLayer} parameterised by $C$.
Hence, stability analysis 
of these equilibria  gives stability of different 
parts of $S(T_a)$.
Specifically, $S(T_a)$ has two normally hyperbolic attracting branches, denoted $S_1(T_a)$ and $S_3(T_a)$, which are separated from a normally hyperbolic repelling branch $S_2(T_a)$ by two non-hyperbolic fold points, denoted $F_1(T_a)$ and $F_2(T_a)$.
The attracting branch $S_1(T_a)$ contains the 
base state $e(T_a)$.
The other attracting branch, $S_3(T_a)$, is the hot metastable state
that arises from the non-monotone microbial respiration function~\eqref{eq:new_r}.\footnote{Note that this stable branch does not exist for the monotone respiration function in~\cite{luke_soil_2011}.}
{\sw 
The aim of GSPT is to combine the slow and fast timescale solutions for $\epsilon=0$, shown in fig.~\ref{fig:frozensystem}~(a), into slow-fast composite solutions for $0<\varepsilon<1$, shown in fig.~\ref{fig:frozensystem}~(b). We point out one particular combination: 
the (blue) {\em special candidate trajectory} $\theta$.
%
}

For $0<\epsilon\ll1$,
the normally hyperbolic parts of $S(T_a)$ are guaranteed 
by the `first' Fenichel theorem
to perturb smoothly to nearby normally hyperbolic and {\em locally invariant  slow manifolds}~\cite{fenichel_1971,fenichel_1979,ckrt_1995,Wechselberger2013,kuehn_2016}. To be specific,  $S_1(T_a)$ and $S_3(T_a)$ perturb to nearby attracting slow manifolds $S_{1,\epsilon}(T_a)$ and $S_{3,\epsilon}(T_a)$. Similarly, $S_{2}(T_a)$ perturbs to a nearby repelling slow manifold $S_{2,\epsilon}(T_a)$.
Each slow manifold is usually non-unique in the sense that there is a family of such slow manifolds that lie exponentially close in $\epsilon$ to each other; see~\cite[Th.3.1]{Wechselberger2013} and~\cite[Th.3.1.4]{kuehn_2016}. The strategy is to fix a representative
for each of these manifolds and work with these representatives.
Near the non-hyperbolic fold points $F_1(T_a)$ and $F_2(T_a)$,
the hyperbolic branches $S_{1,\epsilon}(T_a)$, $S_{3,\epsilon}(T_a)$ and $S_{2,\epsilon}(T_a)$ typically  split, meaning that they typically  become slow manifolds with  either one or two boundaries. 
In particular, $S_{1,\epsilon}(T_a)$ has one inflow boundary near $F_1(T_a)$, $S_{3,\epsilon}(T_a)$ has one outflow boundary near $F_2(T_a)$, while $S_{2,\epsilon}(T_a)$ has one inflow boundary near $F_1(T_a)$ and one outflow boundary near $F_2(T_a)$. 
Local invariance means that trajectories can enter or leave a slow manifold only through its boundary.
Furthermore, the splitting of $S_{1,\epsilon}(T_a)$ and $S_{2,\epsilon}(T_a)$ near $F_1(T_a)$ 
gives rise to
a narrow continuum of special solutions, referred to as {\em canards},
that move along $S_{2,\epsilon}(T_a)$ for some time $t$~\cite{Benoit1983,Dumortier1993,krupaszmolyan2001,SZMOLYAN2001419,Wechselberger2013}{; see fig.~\ref{fig:frozensystem}~(b).}
In the remainder of the paper
we 
follow~\cite{scholarpedia_canard,Benot1981ChasseAC} and~\cite[Sec.3.2.4]{Wechselberger2013}
and identify three different types of canards.
\begin{defn}
\label{defn:canard}
In the autonomous frozen system and in the autonomous compactified system~\eqref{eq: CompactEqns}--\eqref{eq: CompactEqns3}:
\begin{itemize}
    \item[(i)]
    {\em Canards ``without head" } are solutions, or the corresponding trajectories, that move slowly along $S_{2,\epsilon}$ for time $ t= {\cal{O}}(1)$ before moving quickly and directly to $S_{1,\epsilon}$.
    \item[(ii)]
    {\em Canards ``with head"} are solutions, or the corresponding trajectories, that move slowly along $S_{2,\epsilon}$ for time $ t= {\cal{O}}(1)$, then move quickly and directly
    to $S_{3,\epsilon}$, then move slowly along $S_{3,\epsilon}$, before converging to $S_{1,\epsilon}$.
    \item[(iii)] 
    A {\em maximal canard}
    is a special solution, or the corresponding trajectory,  that enters $S_{2,\epsilon}$ through its inflow boundary.
\end{itemize}
\end{defn}
Examples of a maximal canard,  denoted $\theta_\epsilon$, a canard ``without head" and a canard ``with head" in the frozen system are shown in blue, cyan and magenta, respectively, in fig.~\ref{fig:frozensystem}~(b). 
%
{\sw The} (blue) maximal canard {\sw $\theta_\epsilon$} is a special  example of a canard without head, and a perturbation of the special candidate trajectory $\theta$. 
%
In practice, a maximal canard is computed by choosing a suitable arclength and finding the trajectory along $S_2(T_a)$ that takes the longest time to trace out this arclength.

{\sw The above discussion 
identifies two} 
obstacles to defining R-tipping to the hot metastable state. First, the hot metastable state is a transient and thus quantitative phenomenon. In the long term, the  
system converges to the same stable state $e^+$ for any rate $r>0$. Second, there is no unique threshold for R-tipping to the hot metastable state~\cite{Wieczorek2021R}. Rather, one speaks of a ``quasithreshold" comprising a family of  canards.
%
{\sw GSPT} together with compactification allow us to overcome both of these obstacles. 
%

We focus on R-tipping from $e^-$ and relate nonautonomous and compactified system dynamics.
Using the notation $x=(T,C)\in\mathbb{R}^2$ for the state variable of the soil-carbon system, we write
\begin{equation}
\label{eq:pullback}
x^{[r]}(t,e^-)\in\mathbb{R}^2,
\end{equation}
to denote the unique  
solution to the nonautonomous system~\eqref{eq: SmallParam}--\eqref{eq: SmallParam2} at time $t$, that limits to $e^-$ as $t\to -\infty$ for a given rate parameter $r$.\footnote{This solution can be understood as a {\em local pullback attractor} for~\eqref{eq: SmallParam}--\eqref{eq: SmallParam2}~[Th.2.2]\cite{ashwin2017}.}
We also write 
$W^{u,[r]}(\tilde{e}^-)$ 
to denote
the unique rate-dependent one-dimensional unstable invariant manifold of the hyperbolic saddle $\tilde{e}^-$ in the autonomous compactified system~\eqref{eq: CompactEqns}--\eqref{eq: CompactEqns3}.  It follows from~\cite[Prop.6.4(a)]{Wieczorek2021R} 
that  $W^{u,[r]}(\tilde{e}^-)$ contains ${x}^{[r]}(t,{e}^-)$ in the sense that \footnote{For convenience, we leave out the dependence on 
the compactification parameter 
$\nu$ in the notation for the unstable manifold.}
\begin{equation}
\label{eq:Wu}
W^{u,[r]}(\tilde{e}^-)
\supset
\left\{(x,s)~:~ x= x^{[r]}(t,e^-),~s=g_{\nu}(t) \right\}_{t\in\mathbb{R}}.
\end{equation}
This relation allows us to define tracking and R-tipping 
for the nonautonomous system~\eqref{eq: SmallParam}--\eqref{eq: SmallParam2} in terms of 
{\sw $W^{u,[r]}(\tilde{e}^-)$ and
slow manifolds} 
in the autonomous compactified system~\eqref{eq: CompactEqns}--\eqref{eq: CompactEqns3}:
%
\begin{defn}
\label{defn:track}
Consider the nonautonomous system~\eqref{eq: SmallParam}--\eqref{eq: SmallParam2} with exponentially bi-asymptotically constant $T_a(rt)$ and decay coefficient $\rho>0$, and the corresponding autonomous compactified system~\eqref{eq: CompactEqns}--\eqref{eq: CompactEqns3}  with the compactification parameter $\nu\in(0,\rho]$. For a fixed $r>0$:
\begin{itemize}
    \item[(i)] 
    We say $x^{[r]}(t,e^-)$
    {\em tracks} the moving stable equilibrium $e(T_a(rt))$ in system~\eqref{eq: SmallParam}--\eqref{eq: SmallParam2} if $W^{u,[r]}(\tilde{e}^-)$
    connects to $\tilde{e}^+$ directly, 
    i.e.
    without visiting $S_{2,\epsilon}$ or $S_{3,\epsilon}$, 
    in system~\eqref{eq: CompactEqns}--\eqref{eq: CompactEqns3}.
    \item[(ii)] 
    We say $x^{[r]}(t,e^-)$ {\em R-tips} in system~\eqref{eq: SmallParam}--\eqref{eq: SmallParam2}, or say system~\eqref{eq: SmallParam}--\eqref{eq: SmallParam2}
    {\em R-tips} 
    from $e^-$, if $W^{u,[r]}(\tilde{e}^-)$ visits
    $S_{3,\epsilon}$ 
    before 
     connecting to
    $\tilde{e}^+$ in system~\eqref{eq: CompactEqns}--\eqref{eq: CompactEqns3}.
\end{itemize}
\end{defn}

{\sw Since R-tipping is not the complement of tracking\footnote{In the sense that canards "without head" visit $S_{2,\epsilon}$ but not  $S_{3,\epsilon}$ and thus correspond to neither tracking nor R-tipping.},}
the nonautonomous system~\eqref{eq: SmallParam}--\eqref{eq: SmallParam2} with $0<\epsilon\ll 1$ does not have unique R-tipping thresholds separating trajectories that track the moving stable equilibrium from those that R-tip,
and does not have isolated critical rates $r$.
More precisely:
\begin{defn}
\label{defn:critrng}
Consider the nonautonomous system~\eqref{eq: SmallParam}--\eqref{eq: SmallParam2} with exponentially bi-asymptotically constant $T_a(rt)$ and decay coefficient $\rho>0$, and the corresponding autonomous compactified system~\eqref{eq: CompactEqns}--\eqref{eq: CompactEqns3} with the compactification parameter $\nu\in(0,\rho]$.
\begin{itemize}
    \item[(i)]
    We define a {\em critical range of $r$} as an interval of $r$ for which $x^{[r]}(t,e^-)$ neither tracks the moving stable equilibrium {$e(T_a(rt))$} nor R-tips.
    \item[(ii)]
    For a fixed $r>0$, we define an {\em R-tipping quasithreshold}
    in system~\eqref{eq: CompactEqns}--\eqref{eq: CompactEqns3} as a family of canards ``without head'' including a maximal canard, and in system~\eqref{eq: SmallParam}--\eqref{eq: SmallParam2} as a family of solutions $x^{[r]}(t)$ or trajectories corresponding to this family of canards.
     \end{itemize}
\end{defn}
\begin{rmk}
There are important differences between {\em R-tipping quasithresholds} defined above and  {\em regular R-tipping thresholds} introduced in~\cite[Def.5.3]{Wieczorek2021R}.  We refer to section~3 of the electronic appendix for more details.
\end{rmk}
\begin{figure}[t]
    \centering
   \includegraphics[width=0.8\textwidth]{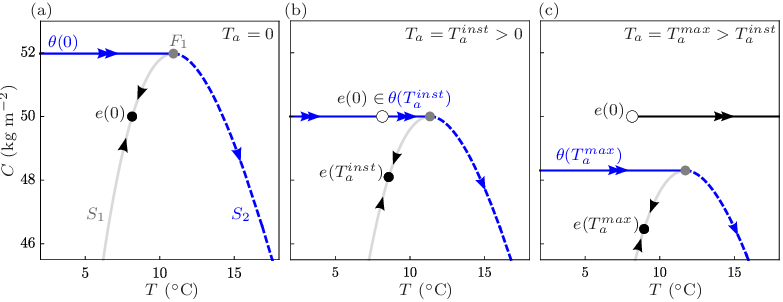}
   \caption{
    Phase portraits of the frozen system with $\epsilon=0$ for three different values of $T_a$ show when to expect R-tipping in the nonautonomous system~\eqref{eq: SmallParam}--\eqref{eq: SmallParam2} with $0 < \epsilon\ll 1$ and a time-varying $T_a$. Note the change in the position of the (blue) special candidate trajectory $\theta(T_a)$ for different values of $T_a$ relative to the fixed initial base state $e(0)$ at $T_a=0$.
   }
    \label{fig:thr_inst}
\end{figure}
%
%
To give intuition for when to expect R-tipping due to crossing a quasithreshold in the nonautonomous system~\eqref{eq: SmallParam}--\eqref{eq: SmallParam2}, we extend the concept
of threshold instability, introduced in~\cite[Def.4.5]{Wieczorek2021R} for regular thresholds, to quasithresholds.
Consider phase portraits of the frozen system 
in the limit
$\epsilon=0$ for three different values of $T_a$ in fig.~\ref{fig:thr_inst}. 
Suppose that $T_a=0^\circ$C and the system is settled at the base state $e(0)$;
see fig.~\ref{fig:thr_inst}~(a).
Then, there is a special value of $T_a$, derived in section 4 of the electronic appendix,
\begin{align}
\label{eq:thr_inst}
T_a = T_a^{inst} \approx -\frac{1}{\alpha}\left(1 + \log{\left(\alpha \frac{A\Pi}{\lambda}\right)} - \alpha \frac{A\Pi}{\lambda}\right),
\end{align}
such that the base state for $T_a=0$ crosses the {\em special candidate trajectory}
for $T_a=T_a^{inst}>0$, that is $e(0)\in \theta(T_a^{inst})$; see fig.~\ref{fig:thr_inst}~(b).
If $T_a$ switches discontinuously from $0$ to $T_a^{max}$, then $e(0)$ becomes the initial condition for the frozen system with $T_a=T_a^{max}$. Thus, {\sw following a switch to }
$T_a^{max}> T_a^{inst}$, $e(0)$ will find itself on the other side of the { \em special candidate trajectory} $\theta(T_a^{max})$ and the system will undergo R-tipping; see fig.~\ref{fig:thr_inst}~(c).  In Sections~\ref{sec:GW} and~\ref{sec:HSA}, we relate condition~\eqref{eq:thr_inst} to R-tipping conditions for $0<\epsilon\ll 1$ and continuously varying $T_a(rt)$.

%

\section{R-tipping Mechanisms for Global Warming}
\label{sec:GW}

\begin{figure}[t]
    \centering
    \includegraphics[width=\textwidth]{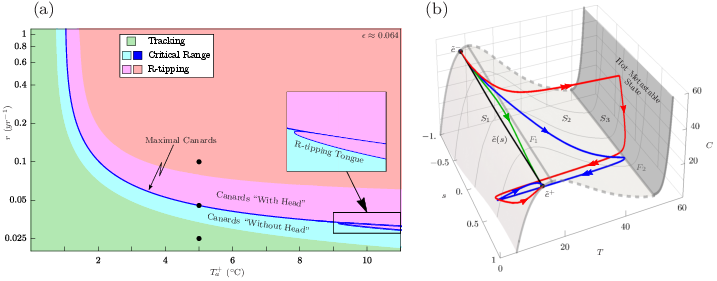}
    \caption{
    (a)
    R-tipping diagram for nonautonomous system~\eqref{eq: SmallParam}--\eqref{eq: SmallParam2} with 
    global warming scenario~\eqref{eq: Tanh Ta}, in the plane of the shift amplitude $T_a^+$ and the rate parameter $r$, for $\epsilon\approx0.064$.
    Shown are regions of (green) tracking, (cyan, blue) critical range, and (magenta, red) R-tipping from $\tilde{e}^{-}$. 
    (b) (Colour) The unstable invariant manifold $W^{u,[r]}(\tilde{e}^-)$ of the saddle $\tilde{e}^-$  for the compactified system~\eqref{eq: CompactEqns}--\eqref{eq: CompactEqns3} with $\epsilon\approx0.064$ and global warming scenario~\eqref{eq: GWCompact} with $T_a^+=5^\circ$C, for three different values of the rate parameter: (green) $r_1=0.02$ , (blue) $r_2\approx0.0454218$, and (red) $r_3=0.1$;
    see the black dots in 
    (a).
    Included for reference is (grey) the critical manifold $S$ defined in~\eqref{eq:CritManGW}. 
    See Table 1 in the electronic appendix for other parameter values.
    }
    \label{fig: GWEps}
\end{figure}

The R-tipping instability in fig.~\ref{fig: Climate Change} arises because time-variation of the  mean atmospheric temperature $T_a$ interacts with the slow timescale of soil carbon $C$. 
{\sw Thus, to uncover the underlying dynamical mechanisms, we consider the 1-fast 2-slow system~\eqref{eq: CompactEqns}--\eqref{eq: CompactEqns3} with the rate parameter $r \lesssim 1$, where $u=rt\lesssim t$ becomes another slow time and $s$ becomes another slow variable.

The singular reduced problem, obtained by setting $\epsilon=0$ for the slow time $t$ in~\eqref{eq: CompactEqns}--\eqref{eq: CompactEqns3},
\begin{align}
        \label{eq: GWrp2}
        \frac{dC}{dt} &= f_2(T,C),\\
        \label{eq: GWrp3}
        \frac{ds}{dt} &= 
        \frac{r}{2}\,
        (1-s^2),
\end{align}
gives slow timescale solutions evolving on a two-dimensional
critical manifold 
\begin{align}
\label{eq:CritManGW}
    S = \left\{(T,C,s)~:~f_1(T,C,T_a^{1}(s)) = 0\right\}\subset \mathbb{R}^2\times[-1,1].
\end{align}
$S$ consists of two normally hyperbolic attracting sheets, $S_1$ containing $\tilde{e}(s)$, and $S_3$ being the hot metastable state, which are separated from a normally hyperbolic repelling sheet $S_2$ by two non-hyperbolic fold curves $F_1$ and $F_2$ (see fig.~\ref{fig: GWEps}~(b)).
The regular layer problem, obtained by setting $\epsilon=0$ for the fast time $\tau = t/\epsilon$,
\begin{align}
\label{eq: GWrp}
\frac{dT}{d\tau} &= f_1(T,C,T_a^{\nu}(s)),
\end{align}
gives fast timescale solutions along straight lines for fixed-in-time $C$ and $s$.

As a model of the global warming scenario, we consider a slow nonlinear shift from $T_a^-=0^\circ$C to a given $T_a^+>0$,
that decays exponentially with the decay coefficient $\rho=2$ as per definition~\ref{defn:expbiconst},
\begin{equation}
\label{eq: Tanh Ta}
    T_a(rt) = \frac{T_a^+}{2}\,(\tanh{(rt)}+1).
\end{equation}
We then fix the compactification parameter $\nu=1$ and apply the inverse of the compactification transformation~\eqref{eq:comptrans} to~\eqref{eq: Tanh Ta}
to obtain\footnote{The inverse of~\eqref{eq:comptrans} is given by equation~(7) in section~2 of the electronic appendix.}
\begin{align}
\label{eq: GWCompact}
T_a^\nu(s) =  
    T_a^{1}(s) =
    \frac{T_a^+}{2}\, \frac{(1 + s)^2}{1 + s^2}.
\end{align}
}

To give an overview of the dynamics near transitions from tracking to R-tipping, we plot in fig.~\ref{fig: GWEps}~(a) an R-tipping diagram in the plane $(T_a^+,r)$
of the input parameters. The diagram was obtained by computing $W^{u,[r]}(\tilde{e}^-)$ in system~\eqref{eq: CompactEqns}--\eqref{eq: CompactEqns3} {\sw with $0< \epsilon \ll 1$} for different values of $T_a^+$ and $r$, and  using
Defs.~\ref{defn:canard}--\ref{defn:critrng} to identify different dynamical regions for system~\eqref{eq: SmallParam}--\eqref{eq: SmallParam2}. 
There are two R-tipping regions, and one tracking-tipping transition found for a large enough shift magnitude $T_a^+$. The smaller R-tipping region is a curious R-tipping tongue.
This tongue gives rise to two (cyan) critical ranges of $r$  for a fixed $T_a^+$ (one of them being very narrow),
and 
different 
{\sw tracking-tipping transitions} 
for low and high magnitudes $T_a^+$ of the shift.

To gain geometric insight into the R-tipping instability caused by
global warming, we plot in fig.~\ref{fig: GWEps}~(b) the unstable manifold $W^{u,[r]}(\tilde{e}^-)$ for a fixed $T_a^+=5^\circ$C and three
different {\sw slow} rates $0<r_1 <r_2<r_3$ (see the black dots in fig.~\ref{fig: GWEps}~(a)), together with the critical manifold {\sw $S$}
for reference. 
%
Figure~\ref{fig: GWEps}~(b) shows that, as $r$ is increased, 
tracking of $\tilde{e}(s)$ by (green) $W^{u,[r_1]}(\tilde{e}^-)$ is lost via canard trajectories, including the maximal canard contained in 
(blue) $W^{u,[r_2]}(\tilde{e}^-)$ that crosses $F_1$ and moves along $S_2$ for the longest time.
This is followed by R-tipping at higher $r$ as shown by (red) $W^{u,[r_3]}(\tilde{e}^-)$ that {\sw crosses $F_1$ and approaches} 
the hot metastable state $S_{3,\epsilon}$ {\sw along the fast $T$-direction} before connecting to $\tilde{e}^+$.
%
{\sw Since R-tipping due to global warming~\eqref{eq: Tanh Ta} occurs when the slow timescale segment of $W^{u,[r]}(\tilde{e}^-)$ 
on $S$ crosses $F_1$,} it should be possible to explain the R-tipping diagram
in fig.~\ref{fig: GWEps}~(a), including the curious R-tipping tongue, in terms of slow 
{\sw timescale solutions of the two-dimensional reduced problem 
alone}~\cite{wieczorek_excitability_2011,perryman_adapting_2014,Vanselow_2019}.

\subsection{
The 
Reduced Problem and Desingularisation} 
\label{sec:GWtm}
The reduced problem{\sw~\eqref{eq: GWrp2}--\eqref{eq: GWrp3}} for the 1-fast 2-slow system~\eqref{eq: CompactEqns}--\eqref{eq: CompactEqns3} with global warming~\eqref{eq: GWCompact} 
describes the evolution of the slow variables $C$ and $s$ in the slow time $t$ on the two-dimensional critical manifold $S$ defined in~\eqref{eq:CritManGW}. However, the onset of R-tipping involves a sudden jump in the fast variable $T$, without any noticeable change in $C$. Thus, it is convenient to consider the evolution of the fast variable $T$ in the slow time $t$ on $S$, which is obtained  by differentiating the 
critical-manifold condition in~\eqref{eq:CritManGW} with respect to $t$. Furthermore,  we use the 
critical-manifold condition  to project the slow flow within $S$ onto the plane $(T,s)$, i.e. eliminate the dependence on $C$, and reformulate the reduced problem as
\begin{align}
\label{eq: ProjRed}
\frac{d T}{dt} &= 
\left.
-
\frac{\partial f_1 /\partial C\cdot f_2 + 
\partial f_1/\partial T_a^1 \cdot dT_a^1 /ds \cdot ds/dt}{\partial f_1/\partial T}
\right|_{S},\\
\label{eq: ProjReds}
\frac{ds}{dt} &= 
 \frac{r}{2}\,
(1-s^2),
\end{align}
where $|_{S}$ denotes restriction to $S$.  
In the reduced problem~\eqref{eq: ProjRed}--\eqref{eq: ProjReds},
the question of loss of tracking boils down to 
whether $W^{u,[r]}(\tilde{e}^-)$ 
leaves $S_1$ through $F_1$.
To address this question, we note that the denominator in~\eqref{eq: ProjRed} changes sign at a fold
\begin{align}
\label{eq:foldsign}
    \left.\frac{\partial f_1}{\partial T}\right|_{S} 
    \left\{
    \begin{array}{ll}
    <0 &\mbox{for}\quad (T,C,s)\in S_1 \cup S_3,\\
    =0 &\mbox{for}\quad (T,C,s)\in F_1 \cup F_2,\\
    >0 &\mbox{for}\quad (T,C,s)\in S_2.
    \end{array}
    \right.
\end{align}
It then follows that, depending on the numerator in~\eqref{eq: ProjRed}, there are two types of fold points:
\begin{itemize}
\item 
A {\em jump point}~\cite{SZMOLYAN2001419} is a  point $p$ on a fold such that
\begin{equation}
\label{eq:jp}
\left.\left(\frac{\partial f_1 }{\partial C}\cdot f_2 + \frac{\partial f_1}{\partial T_a^1} \cdot \frac{dT_a^1}{ds} \cdot \frac{ds}{dt}\right)\right|_{p\in F_1 \cup F_2} \neq 0.
\end{equation}
If a solution of~\eqref{eq: ProjRed}--\eqref{eq: ProjReds} 
approaches a jump point, the denominator in~\eqref{eq: ProjRed} approaches zero while the numerator remains 
finite, so that
$dT/dt$ tends to infinity, and $T(t)$ blows up in $t$. 
This means that the solution reaches a jump point in finite time and ceases to exist within $S$.
Jump points are typically found on open subsets of a fold~\cite{SZMOLYAN2001419}.
\item
A {\em folded singularity}~\cite{SZMOLYAN2001419} is a point $q$ on a fold such that
\begin{equation}
\label{eq:fs}
\left.\left(\frac{\partial f_1 }{\partial C}\cdot f_2 + \frac{\partial f_1}{\partial T_a^\nu} \cdot \frac{dT_a^\nu}{ds} \cdot \frac{ds}{dt}\right)\right|_{q\in F_1\cup F_2} = 0.
\end{equation}
If a solution of~\eqref{eq: ProjRed}--\eqref{eq: ProjReds}
approaches a folded singularity in a special direction where 
the numerator in~\eqref{eq: ProjRed} approaches zero 
faster than or at the same rate as
the denominator in~\eqref{eq: ProjRed},
then $dT/dt$ remains finite, and this solution either grazes or crosses  the fold at $q$ with a finite speed. 
If the crossing is from an attracting to an unstable part of $S$, the solution, or the corresponding trajectory,
is termed a {\em singular canard}~\cite{Benot1981ChasseAC,Benoit1983}.  Folded singularities, which  are typically found as isolated fold points,  are examples of ``singular R-tipping edge states" described in the introduction.
\end{itemize}

Thus, in the singular limit $\epsilon=0$ for the slow time $t$, transitions from tracking to R-tipping caused by
global warming can be understood in terms of singular canards through folded singularities~\cite{wieczorek_excitability_2011,Mitry2013,perryman_adapting_2014}.
The obstacle to analysis of folded singularities and their singular canards is that the right-hand side of~\eqref{eq: ProjRed} is undefined on $F_1$ and $F_2$.
This obstacle is overcome by a {\em desingularisation}~\cite{SZMOLYAN2001419,dumortier1996canard} in the form of a state-dependent time rescaling:
\begin{align}
\label{eq:tresc}
  dt = \left.-d\hat{t}\; \frac{\partial f_1}{\partial T}\right|_{S},
\end{align}
which gives the {\em desingularised system}
%
%
 \begin{align}
\label{eq:DesingDetail}
\frac{d T}{d\hat{t}} &= R_s^*(T)\left(\Pi-\frac{\lambda}{A}\left(T-T_a^{1}(s)\right)\right) + \frac{r\lambda T_a^+}{2A}\left(\frac{1-s^2}{1+s^2}\right)^2,\\
\label{eq:DesingDetail2}
\frac{ds}{d\hat{t}} &= \frac{r\lambda}{2A}(1-s^2)\left(1+ \alpha\left(T-T_a^{1}(s)\right)\frac{ce^{c\alpha(T-b)}-e^{-\alpha(T-b)}}{e^{c\alpha(T- b)}+e^{-\alpha(T- b)}}\right),
\end{align}
defined everywhere on $S$.
The main advantages of 
desingularisation are:
\begin{itemize}
    \item[(a1)] 
    Regular equilibria for the reduced problem~\eqref{eq: ProjRed}--\eqref{eq: ProjReds} remain regular equilibria for the desingularised system~\eqref{eq:DesingDetail}--\eqref{eq:DesingDetail2}.
    Folded singularities for~\eqref{eq: ProjRed}--\eqref{eq: ProjReds} become (new) regular equilibria for~\eqref{eq:DesingDetail}--\eqref{eq:DesingDetail2}. Hence the classification of folded singularities  into ``folded nodes", 
    ``folded foci", 
    ``folded saddles",
    and ``folded saddle-nodes",
    based on their classification into different types of equillibria in~\eqref{eq:DesingDetail}--\eqref{eq:DesingDetail2}~\cite{SZMOLYAN2001419}.
    \item[(a2)]
    According to~\eqref{eq:foldsign} and~\eqref{eq:tresc}, the new time $\hat{t}$ in~~\eqref{eq:DesingDetail}--\eqref{eq:DesingDetail2} flows in the same direction as $t$ on $S_1$ and $S_3$, passes infinitely faster on $F_1$ and $F_2$, and reverses direction on $S_2$. Thus, a phase portrait for~\eqref{eq: ProjRed}--\eqref{eq: ProjReds}
    can be obtained by producing the corresponding phase portrait for~\eqref{eq:DesingDetail}--\eqref{eq:DesingDetail2}, reversing the direction of time (the arrows on trajectories) on $S_{2}$, and relabelling the (new) equilibria on $F_1$ and $F_2$ as folded singularities.
    \item[(a3)]
    It follows from points (a1) and (a2) above that a singular canard through a folded singularity in~\eqref{eq: ProjRed}--\eqref{eq: ProjReds} 
    can be obtained by smoothly concatenating two trajectories tangent to a stable eigendirection of the corresponding equilibrium in~\eqref{eq:DesingDetail}--\eqref{eq:DesingDetail2}.
    We refer to section~5 of the electronic appendix for the discussion of different singular canards associated with different folded singularities.
\end{itemize}

\subsection{Three {\sw Slow} Cases 
in the Reduced Problem}
\label{sec:GW3cases}

To be precise and consistent with Defs.~\ref{defn:track} {\sw and~\ref{defn:critrng}} for $0 < \epsilon \ll 1$, we now define tracking,
R-tipping {\sw and critical rates} for slow external inputs in the limit $\epsilon=0$ as follows:
\begin{defn}
In the reduced problem~\eqref{eq: ProjRed}--\eqref{eq: ProjReds}, {\sw we say that:}
\begin{itemize}
    \item [(i)]
    {\em Tracking} occurs when $W^{u,[r]}(\tilde{e}^-)$ connects to $\tilde{e}^+$ directly, that is without visiting $F_1$.
    \item[(ii)]
    {\em R-tipping} occurs when $W^{u,[r]}(\tilde{e}^-)$ reaches a jump point on $F_1$ {\sw and stops existing in $S$}.\footnote{\sw For $0<\epsilon\ll 1$, the corresponding slow-fast composite solution leaves the attracting slow manifold $S_{1,\epsilon}$ via its outflow boundary and 
    approaches $S_{3,\epsilon}$ along the fast $T$-direction, which is consistent with  Def.~\ref{defn:track} of R-tipping.} 
    \item[(iii)]
    {\sw A {\em critical rate} is an isolated value of $r$ 
    that gives neither tracking nor R-tipping.}
\end{itemize}
\end{defn}
Then, we use relations (a1)--(a3) between the desingularised and reduced system dynamics to identify different 
{\swb \em tracking-tipping transitions}
in the  singular reduced problem~\eqref{eq: ProjRed}--\eqref{eq: ProjReds} through analysis of  the unstable manifold $W^{u,[r]}(\tilde{e}^-)$ in the  regular desingularised system~\eqref{eq:DesingDetail}--\eqref{eq:DesingDetail2}.
We note that, when $r$ is sufficiently small, tracking occurs  because the whole of $F_1$ is repelling, so  $W^{u,[r]}(\tilde{e}^-)$ must remain on $S_1$ and connect directly to $\tilde{e}^+$.
As the rate parameter $r$ is increased, there is a saddle-node bifurcation of equilibria on $F_1$ at some $r=r_{SN}$
in~\eqref{eq:DesingDetail}--\eqref{eq:DesingDetail2}. This bifurcation  corresponds to the appearance of a folded saddle-node type-I singularity (FSN-I)~\cite{KRUPA20102841,VoWechselberger2015} 
in~\eqref{eq: ProjRed}--\eqref{eq: ProjReds}. 
As $r$ is increased past $r_{SN}$, FSN-I bifurcates into a folded saddle (FS) and a folded node (FN). 
%
Thus, $W^{u,[r]}(\tilde{e}^-)$ may interact with singular canards of FS and FN to cross $F_1$ and cause loss of tracking. Analysis of heteroclinic orbits 
in~\eqref{eq:DesingDetail}--\eqref{eq:DesingDetail2}, where $W^{u,[r]}(\tilde{e}^-)$ connects $\tilde{e}^-$ to an equilibrium  on $F_1$, 
reveals three {\swb different} cases of 
{\em loss of tracking}
in~\eqref{eq: ProjRed}--\eqref{eq: ProjReds}:
\begin{itemize}
    \item [(i)] {\em {The} simple {slow} case:}
     $W^{u,[r]}(\tilde{e}^-)$ coalesces with the folded-saddle singular canard $\tilde{\gamma}^S$, and thus crosses from $S_1$ to $S_2$ via FS, at some critical rate $r=r_c>r_{SN}$.
    \item [(ii)] {\em {The} complicated {slow} case:}
     $W^{u,[r]}(\tilde{e}^-)$ 
     grazes $F_1$ at FSN-I when $r=r_{SN}$, 
     coalesces with a  folded-node weak singular canard when $r\in(r_{SN},r_{Ns})$,
     coalesces with the folded-node strong singular canard $\tilde{\gamma}^N_s$ at some $r = r_{Ns}>r_{SN}$,
     and thus crosses from $S_1$ to $S_2$ via FN when $r\in(r_{SN},r_{Ns}]$.
    \item [(iii)] {\em {The} degenerate {slow} case:}
    $W^{u,[r]}(\tilde{e}^-)$ coalesces with the folded-saddle-node singular canard $\tilde{\gamma}^{SN}$,
    and thus crosses from $S_1$ to $S_2$ via FSN-I, at a critical rate $r_c=r_{SN}$.
\end{itemize}
{\sw These three cases of `loss of tracking'  give rise to} {\swb three different cases of `tracking-tipping transition', which we refer to as the {\em three slow cases}.
}
%
Cases (i) and (ii) are typical and closely related to the two cases of non-obvious tipping thresholds described in separate systems in~\cite{perryman_adapting_2014}.
The new case (iii) is special and separates the typical cases (i) and (ii). 
Below, we identify all three cases in an unfolding of a {\em codimension-two non-central heteroclinic FSN-I}, which explains the problem at hand, consolidates the results of~\cite{perryman_adapting_2014}, and is of interest in its own right.
%
%
To 
highlight  interactions of $W^{u,[r]}(\tilde{e}^-)$ with different singular canards, we colour $S_1$
in the phase portraits as follows.
Solutions initialised in the {\em green} region remain in this region and converge directly to $\tilde{e}^+$ {\sw (tracking)}.
Solutions initialised in the {\em dark red}
region reach a jump point of $F_1$ and cease to exist within $S$ {\sw (R-tipping)}.
{\sw Remaining solutions (neither tracking nor R-tipping) include:
weak folded-node singular canards initialised in the {\em yellow} region (singular funnel), the  folded-node strong singular canard $\tilde{\gamma}^{N}_s$ separating the red and yellow regions, and the  folded-saddle singular canard $\tilde{\gamma}^{S}$ separating the green and red regions.
}
%
%
%
Finally, in the figures, we show projections of $S$ onto the plane $(T,s)$.

\subsubsection{The simple {\swb slow} case}

\begin{figure}[t]
    \centering
    \includegraphics[width=8.9cm]{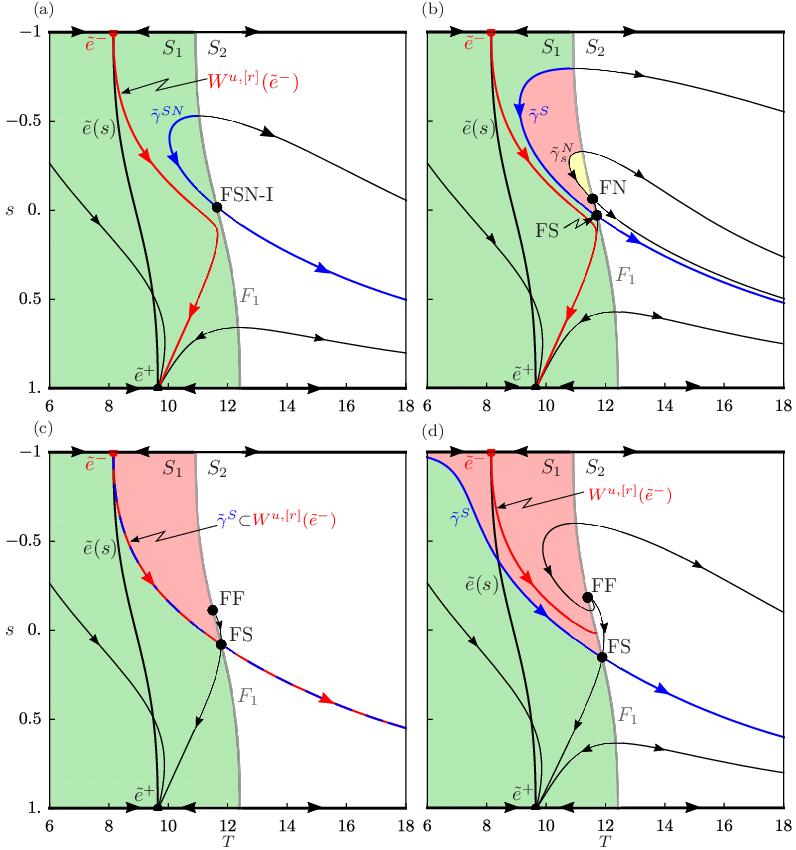}
    \caption{
    {\swb The {\em simple slow 
    case} of tracking-tipping transition in} phase portraits of the reduced problem~\eqref{eq: ProjRed}--\eqref{eq: ProjReds} with nonlinear shift~\eqref{eq: GWCompact}, obtained using the desingularised system~\eqref{eq:DesingDetail}.
    {\swb Note tracking in spite of} a codimension-one local FSN-I in (a), and {\swb loss of tracking via} a codimension-one heteroclinic orbit connecting $\tilde{e}^-$ to FS in (c) {\swb that leads to R-tipping in (d)}.
    $T_a^+=1.5$ and the rate parameter takes values $(a)\ r=r_{SN}\approx0.107194,\ (b)\ r\approx0.108119,\ (c)\ r = r_c\approx0.111459,$ and  $(d)\ r=0.12$; see the black dots in the left inset of fig.~\ref{fig: GWTD}~(c). See Table 1 in the electronic appendix for other parameter values.
    }
    \label{fig:GWSC}
\end{figure}

As $r$ is increased past $r_{SN}$ for $T_a^+=1.5$
in~\eqref{eq: ProjRed}--\eqref{eq: ProjReds}, the appearance of FS and FN via FSN-I gives rise to a folded-saddle singular canard $\tilde{\gamma}^S$, folded-node singular canards including the folded-node strong singular canard $\tilde{\gamma}^N_s$, and new dynamics; {\sw see fig.~\ref{fig:GWSC}~(a)-(b)}. Of particular interest is
an attracting interval of jump points on $F_1$ between FS and FN, and
the corresponding  (dark red) region of 
solutions 
between $\tilde{\gamma}^S$ and $\tilde{\gamma}^N_s$
that reach one of these jump points from $S_1$ and cease to exist within $S$.
%
In this case, $W^{u,[r]}(\tilde{e}^-)$ does not interact with FSN-I, meaning that the FSN-I is local.
$W^{u,[r]}(\tilde{e}^-)$ does not interact with
the folded-node singular canards either because it is separated from them by $\tilde{\gamma}^S$.
Rather, when $r=r_c$, $W^{u,[r]}(\tilde{e}^-)$ and $\tilde{\gamma}^S$ coalesce; see fig.~\ref{fig:GWSC}~(c). At this point, 
tracking is lost 
since $W^{u,[r]}(\tilde{e}^-)$ crosses $F_1$ 
from $S_1$ to $S_2$ via FS. 
For $r>r_c$, $W^{u,[r]}(\tilde{e}^-)$ reaches a jump point on $F_1$ and ceases to exist within $S$
{\sw so R-tipping occurs;}
see fig.~\ref{fig:GWSC}~(d). 

Thus, in the simple {\sw slow} case for $\epsilon=0$, there is 
a {\em critical rate} $r=r_c$. This critical rate is the value of
$r$ that gives a {\it codimension-one hetroclinic orbit} connecting $\tilde{e}^-$ to the corresponding saddle on $F_1$ in the desingularised system~\eqref{eq:DesingDetail}--\eqref{eq:DesingDetail2}.
\subsubsection{The complicated {\swb slow} case}
\begin{figure}[t]
    \centering
   \includegraphics[width=13cm]{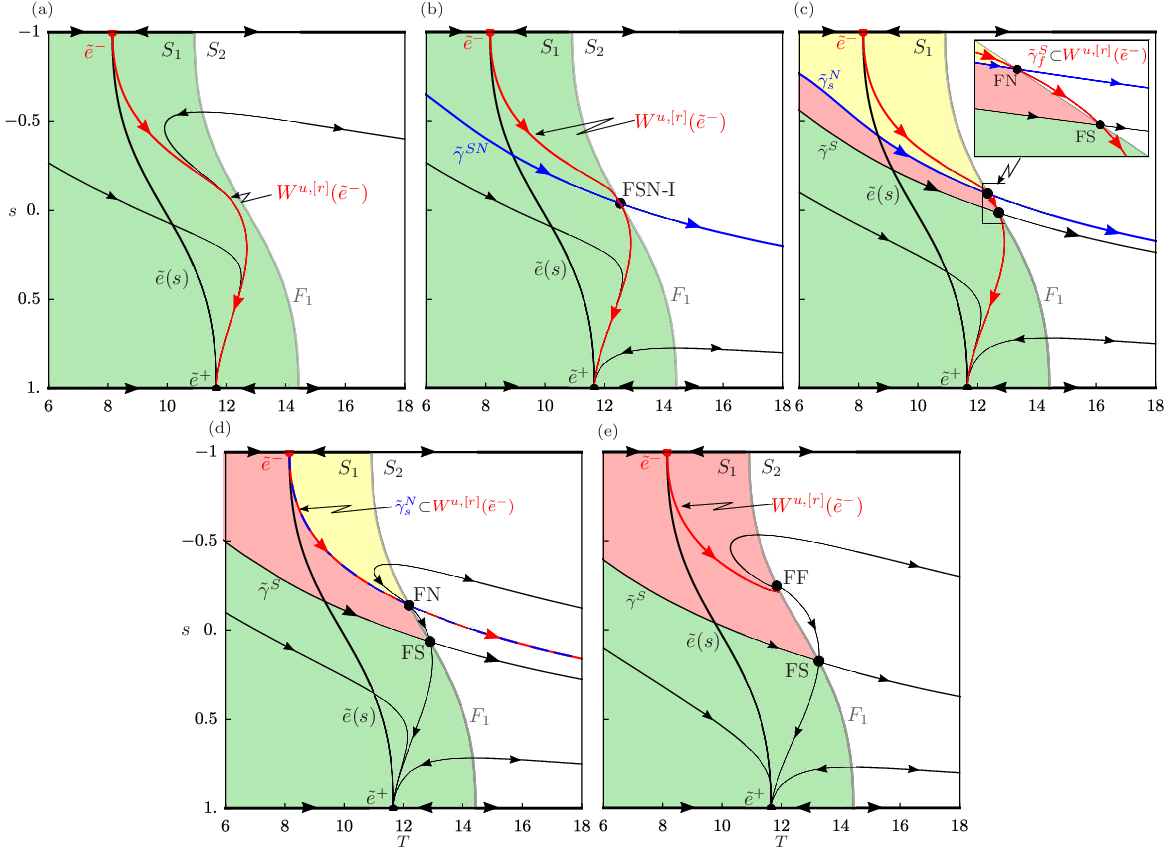}
    \caption{
   {\swb The {\em complicated slow 
   case} of tracking-tipping transition in}
    phase portraits of the reduced problem~\eqref{eq: ProjRed}--\eqref{eq: ProjReds} with nonlinear shift~\eqref{eq: GWCompact}, obtained using the desingularised system~\eqref{eq:DesingDetail}--\eqref{eq:DesingDetail2}.  
    {\swb Note}
    tracking in (a), {\swb loss of tracking via} a codimension-one central heteroclinic FSN-I in (b), and a codimension-one heteroclinic orbit connecting $\tilde{e}^-$ to FN in (d) {\swb that leads to R-tipping in (e)} .
    $T_a^+ = 3.5$ and the rate parameter takes values $(a) \ r = 0.045, (b)\ r = r_{SN} \approx0.050086,\ (c)\ r\approx0.050669,\ (d)\ r = r_c\approx0.052266,$ and  $(e)\ r=0.06$; see the black dots in the right inset of fig.~\ref{fig: GWTD}~(c). See Table 1 in the electronic appendix for other parameter values.
    }
\label{fig:GWCC}
\end{figure}
As $r$ is increased past $r_{SN}$ for $T_a^+=3.5$, the local bifurcation scenario is the same as in the simple 
case{\sw: appearance of FS and FN via FSN-I gives rise to folded-saddle and folded-node singular canards,
an attracting interval of jump points on $F_1$ between FS and FN, and
the corresponding (dark red) region of solutions between $\tilde{\gamma}^S$ and  $\tilde{\gamma}^N_s$ that reach one of these jump points from $S_1$ and cease to exist within $S$; see fig.~\ref{fig:GWCC}. } 
However, the global dynamics are different in that $W^{u,[r]}(\tilde{e}^-)$ interacts with various folded-node singular canards.  
%
%
{\sw When} $r=r_{SN}$, $W^{u,[r]}(\tilde{e}^-)$  
grazes $F_1$ at FSN-I, giving rise to a {\em codimension-one central heteroclinic FSN-I}; see {\eos fig.~\ref{fig:GWCC}~(b)}. 
In other words, the FSN-I is global.
At this point tracking is lost, but R-tipping cannot occur {\sw for $r$ just above $r_{SN}$ because}
%
$W^{u,[r]}(\tilde{e}^-)$
passes through the yellow singular funnel and crosses $F_1$ from $S_1$ to $S_2$ via FN; see {\eos fig.~\ref{fig:GWCC}~(c)}.  Thus, $W^{u,[r]}(\tilde{e}^-)$ contains a  folded-node weak singular canard.
%
Then, $W^{u,[r]}(\tilde{e}^-)$ crosses $F_1$ back to $S_1$ via FS and connects to $\tilde{e}^+$.
Thus, $W^{u,[r]}(\tilde{e}^-)$ also contains the faux folded-saddle singular canard $\tilde{\gamma}^S_f$.
For some 
$r=r_{Ns} > r_{SN}$, $W^{u,[r]}(\tilde{e}^-)$ and $\tilde{\gamma}^N_s$ coalesce. At this point, $W^{u,[r]}(\tilde{e}^-)$ crosses $F_1$ from $S_1$ to $S_2$ via FN 
and never returns to $S_1$; see {\eos fig.~\ref{fig:GWCC}~(d)}.  
It is only for $r>r_{Ns}$ that $W^{u,[r]}(\tilde{e}^-)$ reaches a jump point on $F_1$ and
{\sw R-tipping occurs;}
see {\eos fig.~\ref{fig:GWCC}~(e)}.

Thus, in the  complicated {\sw slow} case for $\epsilon=0$,  there is a {\em critical range of}
$r\in[r_{SN},r_{Ns}]$.
The minimum of the critical range is the value of $r$ that gives a
a {\it codimension-one heteroclinic orbit} connecting $\tilde{e}^-$ to the corresponding saddle-node on $F_1$ along the centre eigendirection in the desingularised system~\eqref{eq:DesingDetail}--\eqref{eq:DesingDetail2}.
The maximum of the critical range is the value of $r$ that gives a
{\it codimension-one hetroclinic orbit} connecting $\tilde{e}^-$ to the corresponding stable node on $F_1$ along the strong eigendirection in~\eqref{eq:DesingDetail}--\eqref{eq:DesingDetail2}.
For $r$ in the interior of the critical range, there is a {\em codimension-zero heteroclinic orbit} connecting $\tilde{e}^-$ to the corresponding stable node on $F_1$ along the weak eigendirection in~\eqref{eq:DesingDetail}--\eqref{eq:DesingDetail2}.

\subsubsection{The degenerate {\swb slow} case}

\begin{figure}[t]
    \centering
    \includegraphics[width=13cm]{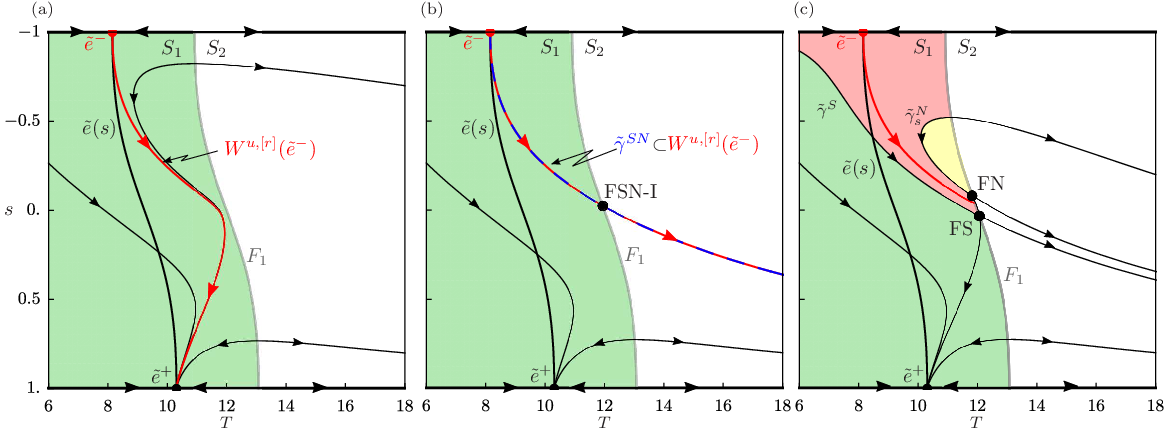}
   \caption{
  {\swb The {\em degenerate slow 
  case} of tracking-tipping transition in}
    phase portraits of the reduced problem~\eqref{eq: ProjRed}--\eqref{eq: ProjReds} with nonlinear shift~\eqref{eq: GWCompact}, obtained using the desingularised system~\eqref{eq:DesingDetail}--\eqref{eq:DesingDetail2}. 
    {\swb Note tracking in (a),} and {\swb loss of tracking via} a codimension-two non-central heteroclinic FSN-I in (b) {\swb that leads to R-tipping in (c)}.  $T_{a}^+=2.15938$ and the rate parameter takes values $(a)\ r\approx0.0716488,\ (b)\ r=r_{SN}=r_c \approx 0.0766488,$ and $(c)\ r\approx0.0776488$; see the black dots in the middle inset of fig.~\ref{fig: GWTD}. System parameter values are given in Table~1 in the electronic appendix. 
    }
    \label{fig:GWDC}
\end{figure}

A natural question to ask is: what 
{\swb separates} the simple and complicated {\sw slow} cases? It turns out that there is a special value $T_a^+=T_{a,c}^+$.
As $r$ is increased past $r_{SN}$ for this special value of $T_a^+$, the local bifurcation scenario is the same as in the simple and complicated {\sw slow} cases. However, the global dynamics are different in that $W^{u,[r]}(\tilde{e}^-)$ interacts with the folded-saddle-node singular canard $\tilde{\gamma}^{SN}$.  
In contrast to the simple and complicated {\sw slow} cases, if $r=r_{SN}$ and $T_a^+=T_{a,c}^+$,  $W^{u,[r]}(\tilde{e}^-)$ and $\tilde{\gamma}^{SN}$ coalesce,
giving rise to a {\em codimension-two non-central heteroclinic FSN-I}; see {\eos fig.~\ref{fig:GWDC}~(b)}. 
In other words, the FSN-I is global, but different from the complicated {\sw slow} case.
At this point, tracking is lost since $W^{u,[r_c]}(\tilde{e}^-)$ crosses $F_1$ from $S_1$ to $S_2$ via FSN-I.
%
For $r>r_{SN}$, $W^{u,[r]}(\tilde{e}^-)$ reaches a jump point on $F_1$ and 
{\sw R-tipping occurs;}
see {\eos fig.~\ref{fig:GWDC}~(c)}.

Thus, in 
the  degenerate case for $\epsilon=0$,  there is 
{\swb neither a critical range nor a critical rate. Instead, there is}
a {\em critical pair} $(T_{a,c}^+,r_c)$.
This critical pair is the combination of $r$ and $T_{a}^+$ that give a
{\em codimension-two heteroclinic orbit} connecting $\tilde{e}^-$ to the corresponding saddle-node on $F_1$ along the stable eigendirection in the desingularised system~\eqref{eq:DesingDetail}--\eqref{eq:DesingDetail2}.\footnote{Note the difference from typical codimension-one heteroclinic orbits connecting $\tilde{e}^-$ to a saddle-node along the centre eigendirection as in the complicated {\sw slow} case.}

\subsection{Bringing it Together: Unfolding of Non-central Heteroclinic FSN-I}
\label{sec:GWbit}

\begin{figure}[t]
    \centering
\includegraphics[width=0.7\textwidth]{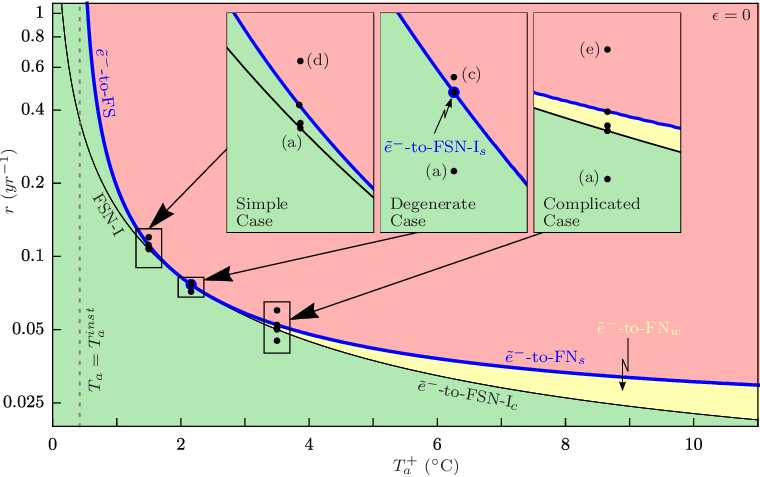}
    \caption{
    Singular R-tipping diagram for the reduced problem~\eqref{eq: ProjRed}--\eqref{eq: ProjReds} with global warming scenario~\eqref{eq: GWCompact}, obtained using the desingularised system~\eqref{eq:DesingDetail}--\eqref{eq:DesingDetail2}. 
    (Green) region of tracking, (red) region of R-tipping, and (yellow) critical range of neither tracking nor R-tipping
    are separated by: (blue) $\tilde{e}^-$-to-FS, $\tilde{e}^-$-to-FSN-I$_s$ and $\tilde{e}^-$-to-FN$_s$, along which $W^{u,[r]}(\tilde{e}^-)$ contains the singular canards $\tilde{\gamma}^S$, $\tilde{\gamma}^{SN}$ and $\tilde{\gamma}^{N}_s$, respectively, and (black) $\tilde{e}^-$-to-FSN-I$_c$, along which $W^{u,[r]}(\tilde{e}^-)$ contains the weak folded saddle-node singular canard. See Table 1 in the electronic appendix for other parameter values.
      }
    \label{fig: GWTD}
\end{figure}

The aim of this section is twofold.
First, we summarise our results from section~\ref{sec:GW}\ref{sec:GW3cases} in the singular R-tipping diagram for $\epsilon=0$ in fig.~\ref{fig: GWTD}~(c). 
Second, we use the singular R-tipping diagram to explain the regular R-tipping diagram for  $\epsilon\approx0.064$ in fig.~\ref{fig: GWEps}~(a).

The singular R-tipping diagram is obtained by the unfolding of a {\em codimension-two non-central heteroclinic
FSN-I} in the plane $(T_a^+,r)$ of the input parameters in the reduced problem~\eqref{eq: ProjRed}--\eqref{eq: ProjReds}. This unfolding, in turn, 
is obtained by the unfolding of a {\em codimension-two non-central saddle-node heteroclinic bifurcation} in the desingularised system~\eqref{eq:DesingDetail}--\eqref{eq:DesingDetail2}.\footnote{This unfolding is reminiscent of the unfolding of a codimension-two non-central saddle-node homoclinic bifurcation in~\cite{chow1990}.}
The  diagram is partitioned into 
{\sw regions of (green) tracking, (red) R-tipping, and (yellow) neither tracking nor R-tipping}
by two curves that are tangent at the special point of  codimension-two non-central hetroclinic FSN-I, denoted $\tilde{e}$-to-FSN-I$_s$ in the middle inset.
Each curve has two branches emanating from the special point.
Both branches of the black curve of FSN-I were obtained by computing the saddle-node bifurcation of equilibria on $F_1$ in~\eqref{eq:DesingDetail}--\eqref{eq:DesingDetail2}; we will return to this curve later.
The left branch of the blue curve,  denoted $\tilde{e}$-to-FS, was obtained by computing a codimension-one heteroclinic orbit connecting $\tilde{e}^-$ to the saddle on $F_1$ in~\eqref{eq:DesingDetail}--\eqref{eq:DesingDetail2}.
This connection gives {\sw a critical rate in} the simple {\sw slow} case for 
the reduced problem~\eqref{eq: ProjRed}--\eqref{eq: ProjReds}, where
$W^{u,[r]}(\tilde{e}^-) \supset \tilde{\gamma}^S$; {\sw see fig.~\ref{fig:GWSC}~(c)}.
%
Note that $\tilde{e}$-to-FS has a vertical asymptote
$T_a^+ = T_a^{inst}\approx 0.423364$, which is given by condition~\eqref{eq:thr_inst}  illustrated in fig.~\ref{fig:thr_inst}.
The right branch of the blue curve, denoted $\tilde{e}$-to-FN$_s$, 
was obtained by computing a codimension-one heteroclinic orbit connecting $\tilde{e}^-$ to the stable node on $F_1$ along the strong eigendirection in~\eqref{eq:DesingDetail}--\eqref{eq:DesingDetail2}. This connection gives the upper boundary of the critical range in the complicated {\sw slow} case for~\eqref{eq: ProjRed}--\eqref{eq: ProjReds}, where $W^{u,[r]}(\tilde{e}^-) \supset \tilde{\gamma}^N_s$; {\sw see fig.~\ref{fig:GWCC}~(d)}.
{\sw
The black curve, denoted $\tilde{e}$-to-FSN-I$_{c}$, was obtained by computing 
a codimension-one heteroclinic orbit connecting $\tilde{e}^-$ to the saddle-node on $F_1$ along the centre eigendirection  in~\eqref{eq:DesingDetail}--\eqref{eq:DesingDetail2}.
This connection gives the lower boundary of the critical range in the complicated {\sw slow} case for~\eqref{eq: ProjRed}--\eqref{eq: ProjReds}, where $W^{u,[r]}(\tilde{e}^-)$ grazes $F_1$ at FSN-I;
{\sw see fig.~\ref{fig:GWCC}~(b)}.\footnote{The left branch of the black curve, denoted FSN-I, does not 
contribute to tracking-tipping transitions.}
}
%
%
%
In the (yellow) region between 
$\tilde{e}$-to-FSN-I$_{c}$ and $\tilde{e}$-to-FN$_s$,  which is the interior of the critical range,
there is a codimension-zero heteroclinic orbit connecting $\tilde{e}^-$ to the stable node on $F_1$ along the weak eigendirection in~\eqref{eq:DesingDetail}--\eqref{eq:DesingDetail2}. 
{\sw Here},
$W^{u,[r]}(\tilde{e}^-)$ contains a  folded-node weak singular canard, meaning that it crosses from $S_1$ to $S_2$ via FN, and also contains
the faux folded-saddle singular canard, meaning that it crosses back to $S_1$ via FS; see the red trajectory in fig.~\ref{fig:GWCC}~(b).
The special point $\tilde{e}$-to-FSN-I$_s$,
where the black and blue curves touch, corresponds to a codimension-two  non-central saddle-node heteroclinic bifurcation that involves a  heteroclinic orbit connecting $\tilde{e}^-$ to the saddle-node on $F_1$ along the stable eigendirection in~\eqref{eq:DesingDetail}--\eqref{eq:DesingDetail2}.
This connection gives
the degenerate case for~\eqref{eq: ProjRed}--\eqref{eq: ProjReds}, where 
$W^{u,[r]}(\tilde{e}^-) \supset \tilde{\gamma}^{SN}$.
Next, we use the singular R-tipping diagram in fig.~\ref{fig: GWTD}~(c) together with rigorous results on 
persistence of singular canards as maximal canards for $0 < \epsilon \ll 1$~\cite{SZMOLYAN2001419,KRUPA20102841,krupaszmolyan2001,Wechselberger2005,VoWechselberger2015} to explain  the regular R-tipping diagram in fig.~\ref{fig: GWEps}~(a). 

In the simple {\sw slow} case, the folded-saddle singular canard $\tilde{\gamma}^S$ perturbs to a family of canards when $0 < \epsilon \ll 1$. The family contains 
one maximal canard, namely
the {\em folded-saddle maximal canard} $\tilde{\gamma}^S_\epsilon$, and  associated canards ``without head" and ``with head"~\cite{SZMOLYAN2001419}.
Thus, 
these canards explain the simple part of the regular R-tipping diagram at lower values of $T_a^+$, including the (cyan) critical range and the vertical asymptote of the blue curve of maximal canards. 
In the context of excitable systems, we use Def.~\ref{defn:critrng}(ii) to relate a family of canards ``without head" associated with a {\em codimension-one normally hyperbolic repelling slow manifold near FS} to a simple
excitability quasithreshold; see also~\cite{wieczorek_excitability_2011,Mitry2013}.

In the complicated {\sw slow} case, the perturbed dynamics for $0 < \epsilon \ll 1$ are far less straightforward 
due to the presence of 
many folded-node maximal canards in addition to the folded-saddle maximal canard.
The number and type of folded-node maximal canards depend on 
the distance between FN and FS, and on
the ratio of the weak and strong eigenvalues of FN, denoted $\mu\in(0,1)$. The folded-node strong singular canard $\tilde{\gamma}^N_s$ perturbs to a family of  canards for all $\mu\in(0,1)$. This family contains the {\em folded-node strong maximal canard} $\tilde{\gamma}^N_{s,\epsilon}$ and associated canards ``without head" and ``with head".
Depending on $\mu$, the family of folded-node weak singular canards from the singular funnel may perturb to a single {\em  folded-node weak maximal canard} $\tilde{\gamma}^{N}_{w,\epsilon}$~\cite{Wechselberger2005,Wechselberger2013}. 
Additionally, there can be a number of {\em secondary folded-node maximal canards} that lie between $\tilde{\gamma}^N_{s,\epsilon}$ and $\tilde{\gamma}^N_{w,\epsilon}$~\cite{SZMOLYAN2001419,Wechselberger2005,Desroches_2010}. 
Furthermore, 
{\sw  the folded-saddle singular faux canard perturbs to folded-saddle faux canards, which may interact with secondary folded-node maximal canards near FSN-I~\cite{Mitry2017}. 
}
The interplay between the folded-node and the folded-saddle can give rise to {\em real-faux canards} described in~\cite{perryman_adapting_2014,VoWechselberger2015}, that are a perturbation of the
red trajectory from fig.~\ref{fig:GWCC}~(b), and to {\em composite canards} identified in~\cite{perryman_adapting_2014}, that follow the folded-saddle maximal canard and then `switch' to and follow one of the folded-node maximal canards.
For example,  we have checked that as $r$ is increased for a fixed $T_a^+ \gtrsim  9.3 ^\circ$C
in fig.~\ref{fig: GWEps}~(a), $W^{u,[r]}(\tilde{e}^-)$ coalesces with a secondary folded-node canard at the bottom boundary of the R-tipping tongue, then with a composite canard at the upper boundary of the tongue, and finally with the folded-node strong maximal canard at 
the bottom boundary of the main R-tipping region{\sw~\cite{osullivan2023}}.\footnote{For higher $T_a^{max}$, we found additional R-tipping tongues, not shown in the figures, whose boundaries involve different secondary and composite canards; see {\sw~\cite[Fig.5.10]{osullivan2023}}.}
Thus, it is folded-node maximal canards and composite canards, together with the associated canards ``without head" and ``with head", that give rise to the complicated tracking-tipping transitions including R-tipping tongue(s), at higher values of $T_a^+$. 
In the context of excitable systems, we use Def.~\ref{defn:critrng}(ii) to relate a family of canards ``without head" associated with a {\em codimension-one normally hyperbolic repelling slow manifold near FN} to a new type of excitability quasithreshold. This quasithreshold is expected to have 
an intricate 
shape as indicated by the computations of slow manifolds  in~\cite{desroches2008, perryman2015}.

\section{The R-Tipping Mechanism for Summer Heatwaves}
\label{sec:HSA}

\begin{figure}[t]
    \centering
    \includegraphics[width=\textwidth]{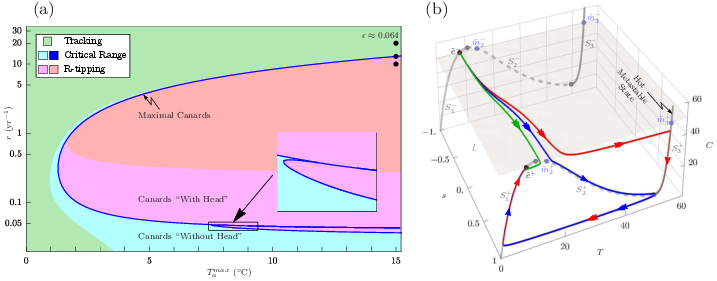}
    \caption{ 
    (a) R-tipping diagram for nonautonomous system~\eqref{eq: SmallParam}--\eqref{eq: SmallParam2} with summer heatwave~\eqref{eq: Sech Ta}, in the plane of the impulse amplitude $T_a^{max}$ and the rate parameter $r$, for $\epsilon\approx0.064${\sw; note the logarithmic scale for $r$.}
    Shown are regions of (green) tracking, (cyan, blue) critical range, and (magenta, red) R-tipping from $\tilde{e}^{-}$. (b)  (Colour) The unstable invariant manifold $W^{u,[r]}_\nu(\tilde{e}^-)$ of the saddle $\tilde{e}^-$ for the compactified system~\eqref{eq: CompactEqns}--\eqref{eq: CompactEqns3} with $\epsilon\approx0.064$, summer heatwave~\eqref{eq:HSACompact} with $T_a^{max}=15^\circ$C and {\eos $\nu = 1/2$}, and three different values of the rate parameter: (green) $r= 20$, (blue) $r\approx12.9123$, and (red) $r=10$; see the black dots in 
    (a). Included for reference are (dark grey) the two disconnected components $S^-$ and $S^+$ of the critical manifold, (light grey) the layer $L$ defined by $C = C^e(0)$, and (light blue) four new equilibria $\tilde{m}_2^\pm$ and $\tilde{m}_3^\pm$ for the layer problem~\eqref{eq: HSALayer}--\eqref{eq: HSALayer2}. See Table 1 in the electronic appendix for other parameter values.
    }
    \label{fig:HSAEps}
\end{figure}

The R-tipping instability in fig.~\ref{fig: Periodic With Hot Summer} arises because time-variation of the atmospheric temperature $T_a$ interacts with the fast timescale of soil temperature $T$. 
%
{\sw Thus, to uncover the underlying dynamical mechanisms, we consider the 2-fast 1-slow system~\eqref{eq: CompactEqns}--\eqref{eq: CompactEqns3} with the rate parameter $r \lesssim 1/\epsilon$, where 
$u= rt \lesssim \tau$
becomes another fast time and $s$ becomes another fast variable.

The singular reduced problem, obtained by setting $\epsilon=1/r =0$ for the slow time $t$ in~\eqref{eq: CompactEqns}--\eqref{eq: CompactEqns3},
\begin{align}
        \label{eq: HSArp}
        \frac{dC}{dt} &= f_2(T,C),
\end{align}
gives slow timescale solutions evolving on a critical manifold that consists of two disconnected one-dimensional components
\begin{align}
\label{eq:CritManHSA}
S = S^- \cup S^+,\quad S^\mp &=\left\{
(T,C,s)~:~f_1(T,C,0)=0,~s=\mp 1
\right\}.
\end{align}
Specifically, 
$
S^- = S^-_1 \cup F_1^- \cup S^-_2 \cup F_2^- \cup S^-_3 
$ 
is the critical manifold of the past limit system~\eqref{eq:ls-} with $T_a^-=0^\circ$C, embedded in the compactified phase space where it gains one unstable direction. Similarly, 
$
S^+ = S^+_1 \cup F_1^+ \cup S^+_2 \cup F_2^+ \cup S^+_3
$ 
is the critical manifold of the future limit system~\eqref{eq:ls+} with $T_a^+=0^\circ$C, embedded in the compactified phase space where it gains one stable direction. 
The attracting branch $S_3^+$ 
is the hot metastable state;
see fig.~\ref{fig:HSAEps}~(b).
%
The 
layer problem, obtained by
{\swr defining $r = \tilde{r}/\epsilon$ for some 
constant $\tilde{r}={\cal O}(1)$ and taking
the limit 
$\epsilon\to 0$}
for the fast time $\tau = t/\epsilon$,
\begin{align}
\label{eq: HSALayer}
        \frac{dT}{d\tau} &= f_1(T,C,T_a^{\nu}(s)),\\
\label{eq: HSALayer2}
        \frac{ds}{d\tau} &=  
{\swr \frac{\tilde{r}\,\nu}{2}}\,(1-s^2),
\end{align}
gives fast timescale solutions on a two-dimensional {\em layer} with a fixed-in-time $C$,
\begin{equation}
\label{eq:L}
L=\left\{(T,C, s)~:~C = {\rm const.} \right\}.  
\end{equation}

As a model of a summer heatwave, we sacrifice seasonal variations to simplify the analysis\footnote{An extension of the compactification to asymptotically periodic inputs, such as in fig.~\ref{fig: Periodic With Hot Summer}(a), is left for future research.}
and consider a fast impulse rising from $T_a^{min}=0^\circ$C to a given $T_a^{max}>0^\circ$C and then dropping back to $0^\circ$C, that decays exponentially with the decay coefficient $\rho=1$ as per definition~\ref{defn:expbiconst},
\begin{equation}
\label{eq: Sech Ta}
    T_a(rt) = T_a^{max}\sech{(rt)}.
\end{equation}
We then fix the compactification parameter $\nu= 1/2$ and apply the inverse of the compactification transformation~\eqref{eq:comptrans} to~\eqref{eq: Sech Ta} to obtain\footnote{The inverse of~\eqref{eq:comptrans} is given by equation~(7) in section~2 of the electronic appendix.}
\begin{align}
\label{eq:HSACompact}
   T_a^\nu (s) =T_a^{\frac{1}{2}}(s) = \frac{2\,T_a^{max}\,(1-s)^{2}}{(1+s)^{4}+(1-s)^{4}}.
\end{align}
}

To give a full overview of 
transitions from tracking to R-tipping for impulse inputs, and to make connections to section~\ref{sec:GW}, we plot an R-tipping diagram in the plane $(T_a^{max},r)$ 
of the input parameters in fig.~\ref{fig:HSAEps}~(a) for a wide range of the rate parameter $r$. 
The diagram was obtained by computing $W^{u,[r]}(\tilde{e}^-)$ in system~\eqref{eq: CompactEqns}--\eqref{eq: CompactEqns3} {\sw with $0< \epsilon \ll 1$} for different values of $T_a^{max}$ and $r$, and using Defs.~\ref{defn:canard}--\ref{defn:critrng} to identify different dynamical regions for system~\eqref{eq: SmallParam}--\eqref{eq: SmallParam2}. 
{\sw The families of canards were considered relative to the evolving slow manifold, from  two-dimensional $S$ in~\eqref{eq:CritManGW} to one dimensional $S^+$ in~\eqref{eq:CritManHSA}, as $r$ was increased.} 
{\sw The} shape of the R-tipping region in fig.~\ref{fig:HSAEps}~(a) is rather different to that obtained in fig.~\ref{fig: GWEps}~(a). 
There are two R-tipping tongues, one large akin to that in~\cite[Sec.6]{okeeffe_tipping_2020} and one small akin to that in fig.~\ref{fig: GWEps}~(a), each enclosing a separate (magenta, red) region of R-tipping.
%
%
The lower part of the diagram corresponds to slow impulses that last for decades {\sw (1-fast 2-slow system)}.
The tracking-tipping transitions 
in this part occur during the rise of the impulse $T_a(rt)$, and 
{\sw are} of the same type as {\sw those}
described in section~\ref{sec:GW}.
The upper part of the diagram corresponds 
to fast impulses {\sw (2-fast 1-slow system)} and is the focus of this section. 
{\sw The} rate parameter $r$ in the range of 10 to 15 yr$^{-1}$ gives a pulse duration in the range of 3 to 2 months, in line with realistic summer heatwaves~\cite{perkins-kirkpatrick_lewis_2020}.\footnote{The pulse duration is obtained using the formula 
$2\ln(2 + \sqrt{3})/r \approx 2.6/r$ for the full width at half maximum of the hyperbolic secant $\sech(rt)$.} Thus, the single tipping-tracking transition found in the upper part of the diagram quantifies the intensity and duration of  summer heatwaves that trigger R-tipping to the hot metastable state in the soil-carbon system~\eqref{eq: SmallParam}--\eqref{eq: SmallParam2}.

To gain geometric insight into the R-tipping instability caused by
a summer heatwave, we 
plot in fig.~\ref{fig:HSAEps}~(b) the unstable manifold $W^{u,[r]}(\tilde{e}^-)$ for a fixed $T_a^+=15^\circ$C and three different {\sw fast} rates $r_1 >r_2>r_3 >0$ (see the black dots 
in fig.~\ref{fig:HSAEps}~(a)), {\sw 
together with the (light gray) layer $L$ of constant $C = C^e(0)$ and the (dark gray) critical manifold $S = S^-\cup S^+$.
}
%
%
%
Figure~\ref{fig:HSAEps}~(b) shows that,
as $r$ is {\em decreased} (the duration of the heatwave is {\em increased}), 
tracking by (green) $W^{u,[r_1]}(\tilde{e}^-)$ is lost via canard trajectories, including the maximal canard contained in (blue) $W^{u,[r_2]}(\tilde{e}^-)$ that  approaches
$S_2^+$ 
and moves along $S_2^+$ for the longest time. 
This is followed by R-tipping at {\em lower} $r$ as shown by (red) $W^{u,[r_3]}(\tilde{e}^-)$ that visits the hot metastable state $S^+_{3,\epsilon}$ before connecting to $\tilde{e}^+\in S^+_{1,\epsilon}$.
{\sw Since tracking-tipping transitions due to a summer heatwave~\eqref{eq: Sech Ta} occur when the fast timescale segment of $W^{u,[r]}(\tilde{e}^-)$ 
on $L$ approaches the hyperbolic saddle branch $S_2^+$ of $S^+$,} 
it should be possible to explain {\sw the upper part of} the R-tipping 
diagram in fig.~\ref{fig:HSAEps} in terms of 
{\sw fast timescale solutions of a suitably chosen two-dimensional layer problem alone.}


     %

\subsection{One {\sw Fast} 
Case in the Layer Problem 
}
\label{sec:HSAtlp}

The suitable layer problem for the 2-fast 1-slow system~\eqref{eq: CompactEqns}--\eqref{eq: CompactEqns3} with a summer heatwave~\eqref{eq:HSACompact} is obtained by {\sw fixing $C$ in~\eqref{eq: HSALayer}--\eqref{eq: HSALayer2} and~\eqref{eq:L} at the equilibrium soil-carbon concentration $C^e(0)$ for the past and future limit systems. 
}
%
Crucially, 
{\sw such a}
layer problem~\eqref{eq: HSALayer}--\eqref{eq: HSALayer2} has six equilibria.
Two of these equilibria, namely the saddle $\tilde{e}^- \in S^-_1$ and the {\em {\sw base-state} sink} $\tilde{e}^+ \in S^+_1$, are the equilibria of the compactified system~\eqref{eq: CompactEqns}--\eqref{eq: CompactEqns3}. 
The four new equilibria, that do not exist in~\eqref{eq: CompactEqns}--\eqref{eq: CompactEqns3}, are found at the intersections of $L$ with $S^\pm_{2,3}$. 
These are: 
the source $\tilde{m}^-_2\in S_2^-$,
the saddle $\tilde{m}^-_3\in S_3^-$,
the saddle $\tilde{m}^+_2\in S_2^+$, and 
the {\em {\sw hot-state} sink} $\tilde{m}^+_3\in S_3^+$; see the (blue) dots in fig.~\ref{fig:HSAEps}~(b).

\begin{figure}[t]
\centering
    \includegraphics[width=0.9\textwidth]{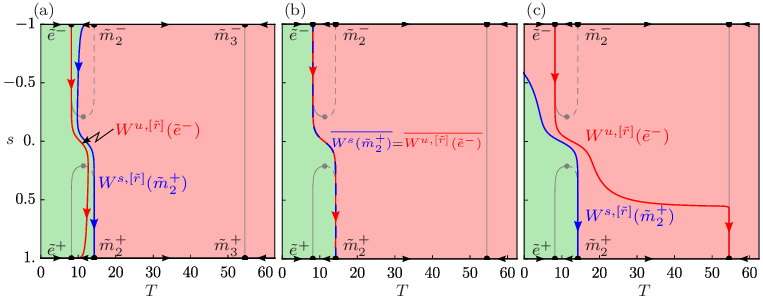}
\caption{
{\swb The {\em fast 
case} of tracking-tipping transition in} phase portraits of the layer problem~\eqref{eq: HSALayer}--\eqref{eq: HSALayer2} {\sw with $C=C^e(0)$} 
{\eos and $\nu=1/10$}.
{\swb Note tracking in (a), and loss of tracking via} a codimension-one heteroclinic orbit connecting $\tilde{e}^-$ to $\tilde{m}_2^+$ in (b) {\swb that leads to R-tipping in (c)}. 
$T_a^+=15$ and the rate parameter takes values 
{\eos $(a)\ \tilde{r} = 1.28$, $(b)\ \tilde{r}=\tilde{r}_c\approx0.972269$, and $(c)\ \tilde{r} = 0.64$;}
see the black dots in the inset of fig.~\ref{fig:HSAtd}. 
Included for reference are different branches of (grey) $\tilde{M}$.
 See the electronic appendix for  other parameter values in Table 1 and  for definition of $\tilde{M}$ in section 6.
}
\label{fig:HSAfc}
\end{figure}

To be precise and consistent with
Defs.~\ref{defn:track} {\sw and~\ref{defn:critrng}} for $0 < \epsilon \ll 1$, we now define tracking,
R-tipping {\sw and critical rates} 
for fast external inputs in the limit $\epsilon=0$ as follows:
\begin{defn}
In the layer problem~\eqref{eq: HSALayer}--\eqref{eq: HSALayer2}, 
{\sw we say that},
\begin{itemize}
    \item [(i)]
    {\em Tracking} occurs when $W^{u,[{\eos \tilde{r}}]}(\tilde{e}^-)$ connects to {\sw the base-state sink} $\tilde{e}^+$.
    \item[(ii)]
    {\em R-tipping} occurs when $W^{u,[{\eos \tilde{r}}]}(\tilde{e}^-)$ connects to {\sw the hot-state sink} $\tilde{m}_3^+$.
    \item[(iii)]
    {\sw A {\em critical rate} is an isolated value of ${\eos \tilde{r}}$ 
    that gives neither tracking nor R-tipping.
    }
\end{itemize}
\end{defn}
In other words, in the the layer problem~\eqref{eq: HSALayer}--\eqref{eq: HSALayer2}, the question of loss of tracking boils down to whether or not $W^{u,[{\eos\tilde{r}}]}(\tilde{e}^-)$ connects to the saddle $\tilde{m}_2^+$.
Thus, the saddle $\tilde{m}_2^+$ is another example of the ``singular R-tipping edge state" described in the introduction. 

The  dynamics of the layer problem~\eqref{eq: HSALayer}--\eqref{eq: HSALayer2} on $L$
are shown  in a series of phase portraits in fig.~\ref{fig:HSAfc}, where we fix $T_a^{max} = 15 > T_a^{inst}$ and decrease 
{\eos $\tilde{r}$} 
across the critical rate.
The basin of attraction of $\tilde{e}^+$ is plotted in {\em green}, and the basin of attraction of $\tilde{m}_3^+$ is plotted in {\em dark red}.
{\sw The rate-dependent (blue) stable manifold of $\tilde{m}_2^+$, denoted} $W^{s,[{\eos \tilde{r}}]}(\tilde{m}_2^+)$, is contained in the basin boundary separating these two basins of attraction.
For ${\eos \tilde{r}}$ sufficiently large, (red) $W^{u,[{\eos \tilde{r}}]}(\tilde{e}^-)$ lies in the basin of attraction of $\tilde{e}^+$
and connects to $\tilde{e}^+$; see fig.~\ref{fig:HSAfc}~(a). 
When ${\eos \tilde{r}=\tilde{r}_c}$, $W^{u,[{\eos \tilde{r}}]}(\tilde{e}^-)$ coalesces with $W^{s,[{\eos \tilde{r}}]}(\tilde{m}_2^+)$ along the basin boundary. At this point, tracking is lost since $W^{u,[{\eos \tilde{r}}]}(\tilde{e}^-)$ connects to the saddle $\tilde{m}_2^+$; see fig.~\ref{fig:HSAfc}~(b). 
For ${\eos \tilde{r}<\tilde{r}_c}$, $W^{u,[{\eos \tilde{r}}]}(\tilde{e}^-)$ lies in the basin of attraction of $\tilde{m}_3^+$
and connects to $\tilde{m}_3^+$.
In other words, 
{\sw R-tipping occurs}
for $r<r_c$; see fig.~\ref{fig:HSAfc}~(c). 

Thus, in the layer problem~\eqref{eq: HSALayer}--\eqref{eq: HSALayer2} {\sw with $C=C^e(0)$}, 
{\sw there is just one case of a tracking-tipping transition, where tracking is lost at a}
{\em critical rate} ${\eos \tilde{r}=\tilde{r}}_c$. 
{\sw We refer to this case as the {\swb \em fast case} to distinguish it from the three {\em slow cases} identified in section~\ref{sec:GW}\ref{sec:GW3cases}.}
The critical rate is the value of $r$ that gives a {\it codimension-one hetroclinic orbit} connecting $\tilde{e}^-$ to the new saddle $\tilde{m}_2^+$ for the layer problem~\eqref{eq: HSALayer}--\eqref{eq: HSALayer2}. 
%
For additional insight into R-tipping caused by a summer heatwave, we refer to section~6 of the electronic appendix.

\subsection{Bringing it Together: Additional Saddle-to-Saddle Heteroclinic Orbit}
\label{sec:HSAbiat}

\begin{figure}[t]
    \centering
    \includegraphics[width=0.8\textwidth]{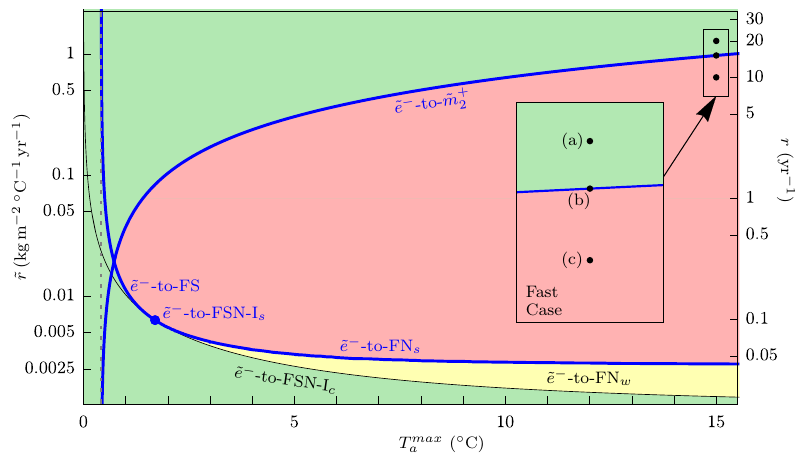}
    \caption{
    Singular R-tipping diagram obtained using a combination of {\sw two problems}.
    (The upper curve $\tilde{e}^-$-to-$\tilde{m}^+_2$) the {\sw two-dimensional} layer problem~\eqref{eq: HSALayer}--\eqref{eq: HSALayer2} with summer heatwave~\eqref{eq:HSACompact}, {\swr rate parameter $\tilde{r}=\epsilon r$ shown on the left vertical axis, and $\nu=1/2$}.
    (The lower curves) the {\sw two-dimensional} reduced problem~\eqref{eq: ProjRed}--\eqref{eq: ProjReds} with slow impulse~\eqref{eq:HSACompact} {\swr and rate paraemeter $r$ shown on the right vertical axis}. 
    Note the logarithmic {\swr scales} for {\eos $r$ and $\tilde{r}$.}
    At higher ${\eos \tilde{r}}$, the (red) region of R-tipping is separated from the (green) region of tracking by {\swr the} (blue) {\swr curve} $\tilde{e}^-$-to-$\tilde{m}^+_2$, along which
    $W^{u,[{\eos \tilde{r}}]}(\tilde{e}^-)$ intersects $W^{s,[{\eos \tilde{r}}]}(\tilde{m}^+_2)$ in~\eqref{eq: HSALayer}--\eqref{eq: HSALayer2}.
    For lower $r$, see the caption of fig.~\ref{fig: GWTD}~(c).
    See Table 1 in the electronic appendix for other parameter values.
    }
    \label{fig:HSAtd}
\end{figure}

%
{\eos In this section, we first}
produce the singular R-tipping diagram for $\epsilon = 0$ in fig.~\ref{fig:HSAtd},
{\swr and then use it}
to explain the regular R-tipping diagram for $\epsilon\approx0.064$ in fig.~\ref{fig:HSAEps}~(a).

In the singular R-tipping diagram {\sw in fig.~\ref{fig:HSAtd}}, the lower-$r$  tracking-tipping transitions were obtained in the same manner as outlined in section~\ref{sec:GW}, that is by computing 
different heteroclinic orbits in 
{\sw the desingularised system for the slow impulse $T_a^{1/2}(s)$ in~\eqref{eq:HSACompact}.}
To avoid repetition,  we skip the derivations and refer to sections~\ref{sec:GW}\ref{sec:GW3cases} and~\ref{sec:GW}\ref{sec:GWbit} for more details.

The higher-${\eos \tilde{r}}$ tipping-tracking transition, which is the focus of this section, is given by the (blue) curve denoted $\tilde{e}^-$-to-$\tilde{m}_2^+$. This curve was obtained by computing a 
codimension-one heteroclinic orbit connecting the saddle $\tilde{e}^-$ to the saddle $\tilde{m}_2^+$ in the layer problem~\eqref{eq: HSALayer}--\eqref{eq: HSALayer2}; see fig.~\ref{fig:HSAfc}~(b).  This {\sw heteroclinic} connection 
{\sw gives the critical rate when $\overline{W^{u,[{\eos \tilde{r}}]}(\tilde{e}^-)} = \overline{W^{s,[{\eos \tilde{r}}]}(\tilde{m}_2^+)}$ \footnote{We write $\overline{A}$ to denote the closure of $A$, that is the smallest closed set containing $A$.}, and is} the only case 
{\sw of loss of tracking} for~\eqref{eq: HSALayer}--\eqref{eq: HSALayer2}.
%
%
Note that curves $\tilde{e}^-$-to-$\tilde{m}_2^+$ and $\tilde{e}^-$-to-FS have a common vertical asymptote $T_a^{max} = T_a^{inst} \approx 0.423364$, which is given by condition~\eqref{eq:thr_inst}  illustrated in fig.~\ref{fig:thr_inst}.

Next, in addition to the `first' Fenichel theorem on persistence of normally hyperbolic critical manifolds as slow manifolds for $0 < \epsilon \ll 1$, we recall the `second' Fenichel theorem that guarantees persistence of stable and unstable invariant manifolds of normally hyperbolic saddle critical manifolds when $0 < \epsilon \ll 1$~\cite{fenichel_1971,fenichel_1979,ckrt_1995,Wechselberger2013,kuehn_2016}.
We then use the singular R-tipping diagram in fig.~\ref{fig:HSAtd}, together with both Fenichel theorems and rigorous results on 
persistence of singular canards as maximal canards for $0 < \epsilon \ll 1$~\cite{SZMOLYAN2001419,KRUPA20102841,krupaszmolyan2001,VoWechselberger2015,Wechselberger2005}, to explain  the regular R-tipping diagram in fig.~\ref{fig:HSAEps}~(a). 
%

In the {\em fast case}, the heteroclinic $\tilde{e}^-$-to-$\tilde{m}_2^+$ connection for the layer problem~\eqref{eq: HSALayer}--\eqref{eq: HSALayer2} perturbs to an intersection of  
$W^{u,[{\eos \tilde{r}}]}(\tilde{e}^-)$ with the (fast) stable manifold of the normally hyperbolic {\em saddle slow manifold} $S_{2,\epsilon}^+$ in the compactified system~\eqref{eq: CompactEqns}, giving rise to the (blue) maximal canard shown in fig.~\ref{fig:HSAEps}~(b).  In other words, 
$W^{u,[{\eos \tilde{r}}]}(\tilde{e}^-)$ approaches the normally hyperbolic saddle slow manifold
$S^+_{2,\epsilon}$ along its stable manifold, thus moves along $S^+_{2,\epsilon}$ for as long $S^+_{2,\epsilon}$ exists, then approaches the stable slow manifold $S^+_{1,\epsilon}$
and eventually connects to $\tilde{e}^+$. In addition, for nearby values of $r$, we 
expect associated canards ``without head'' and canards``with head''. These occur when
$W^{u,[{\eos \tilde{r}}]}(\tilde{e}^-)$ near misses the stable manifold of  $S_{2,\epsilon}^+$ to the left or right, respectively, and diverges away from $S_{2,\epsilon}^+$ before $S_{2,\epsilon}^+$ ceases to exist.
Thus, it is the (blue) maximal canard shown in fig.~\ref{fig:HSAEps}~(b), together with the associated canards ``without head",  that give rise to the higher-$r$ tipping-tracking transition
in the regular R-tipping diagram in fig.~\ref{fig:HSAEps}~(a).
In the context of excitable systems, we use Def.~\ref{defn:critrng}(ii) to relate a family of canards ``without head" associated with a {\em codimension-one stable manifold of a normally hyperbolic saddle slow manifold} to yet another new type of excitability  quasithreshold.

\section{Conclusion and Outlook}
In this paper, we obtain two kinds of results.
First, we demonstrate that sufficiently fast atmospheric warming can cause  R-tipping to a subsurface hot metastable state in bioactive peat soils, using a conceptual process-based ODE model with realistic soil parameter values and contemporary climate patterns. This gives rise to the hypothesis that such R-tipping is a main cause of ``Zombie fires" observed in peatlands, that disappear from the surface, smoulder underground during the winter, and ``come back to life" in the spring. 
Second, we recognise that such R-tipping is a nonautonomous instability that occurs due to crossing an elusive quasithreshold in the phase space of a multiple-timescale dynamical system, and thus cannot be explained by traditional autonomous stability theory.  Therefore, to understand the underlying dynamical mechanisms, we provide a mathematical framework that is underpinned by a compactification technique for asymptotically autonomous dynamical systems and concepts from geometric singular perturbation theory, such as slow manifolds, canard trajectories and folded singularities.
This {\sw framework} explains R-tipping to the hot metastable state in the soil-carbon system and, more generally, identifies generic cases of  {\sw tracking-tipping transitions} 
due to crossing a quasithreshold. Furthermore, {\sw we show} 
that a quasithreshold  gives rise to critical ranges of the rate of change of the external input rather than isolated critical rates, and reveal new types of quasithresholds. These results pose two types of challenges for future research.

First, to strengthen our hypothesis, it would be interesting to build on the excellent agreement with the medium-complexity PDE model in fig.~\ref{fig: Khvorostyanov Comparison} and extend the conceptual ODE model in different directions. For example, 
    include terms describing peat soil ignition processes at higher temperatures~\cite{yuan_restuccia_rein_2021}. This would allow investigation of existence of additional (quasi)thresholds for the onset of flames.
    Include other physical processes, in addition to soil temperature, that contribute to microbial soil respiration. For example, 
    the inclusion of soil hydrology would extend the applicability of the model; see for example DigiBog\_Hydro~\cite{Digibog_hydro_2022} and~\cite{Waddington2014}.
    Include random fluctuations to account for the fact that climate and weather patterns are inherently noisy. This would give a more accurate description of the tipping being a combination of R-tipping and noise-induced tipping (N-tipping). Understanding how the two tipping mechanisms interact is non-trivial~\cite{Ritchie2017,Slyman2022} and an area of ongoing research for  quasithresholds.
    Account for spatio-temporal dynamics by extending the conceptual ODE model to a PDE model. Vertical diffusion has already been studied in the extended Luke and Cox model with interesting results confirming validity of the conceptual ODE model~\cite{Clarke2021}. 
    Finally, it might be interesting to couple the conceptual model to spatially extended land-surface models such as JULES~\cite{JULES_2022} or DigiBog~\cite{Digibog_2022}. Currently, such models often neglect heat produced by microbial decomposition which is a key factor for the R-tipping instability to the hot metastable state described in this paper.

Second, there are challenges related to extending the mathematical framework and obtaining rigorous results.
{\sw One is {\eos the} extension of the compactification technique to asymptotically periodic inputs and usage of geometrical singular perturbation theory to analyse external inputs with seasonal variations.}
{\sw Another is a general theory of R-tipping due to crossing quasithresholds in multiple-timescale systems, for which} the definitions of R-tipping and a quasithreshold introduced here for the soil-carbon system are a good starting point. 
Of particular interest to scientists and mathematical modelers are rigorous yet easily testable criteria for the occurrence of such R-tipping, akin to those derived in~\cite{Wieczorek2021R} for R-tipping due to crossing regular thresholds.

\renewcommand{\abstractname}{Competing Interests}
\begin{abstract}
 We declare we have no competing interests.
\end{abstract}

\renewcommand{\abstractname}{Funding}
\begin{abstract}
EO'S acknowledges financial support from ATSR Ltd.
SW acknowledges partial support by the EvoGamesPlus Innovative Training Network funded by the European Union’s Horizon 2020 research and innovation programme under the Marie Skłodowska-Curie grant agreement No 955708.
KM acknowledges partial support from Enterprise Ireland grants IP 2018 0747 and IP 2021 1008.
\end{abstract}

\renewcommand{\abstractname}{Acknowledgments}
\begin{abstract}
The authors would like to thank Hassan Alkhayoun, {\sw Martin Wechselberger and two anonymous referees} for 
{\sw their constructive}
comments. 
\end{abstract}



\bibliographystyle{RS2}
\bibliography{main}

\end{document}